\documentclass []{article}
\usepackage{amssymb}
\def\qmod#1#2{{\hbox{}^{\displaystyle{#1}}}\!\big/\!\hbox{}_{
\displaystyle{#2}}}

\def\resto#1#2{{
#1\hskip 0.4ex\vline_{\hskip 0.4ex\raisebox{-1ex}
{{${\scriptstyle #2}$}}}}}

\newfam\msbfam

\font\tenmsb=msbm10

\font\sevenmsb=msbm10 at 7pt
\font\fivemsb=msbm10 at 5pt

\textfont\msbfam=\tenmsb
\scriptfont\msbfam=\sevenmsb
\scriptscriptfont\msbfam=\fivemsb

\def\Bbb{\fam\msbfam\tenmsb}

\def\C{{\Bbb C}}

\def\H{{\Bbb H}}

\def\P{{\Bbb P}}
\def\Q{{\Bbb Q}}
\def\R{{\Bbb R}}
\def\T{{\Bbb T}}
\def\Z{{\Bbb Z}}

\def\cringle{\mathaccent23}
\def\union{\mathop{\bigcup}}
\def\qed {\hfill\vrule height6pt width6pt depth0pt \bigskip}

\def\map{\longrightarrow}
\def\textmap#1{\mathop{\vbox{\ialign{
                                 ##\crcr
     ${\scriptstyle\hfil\;\;#1\;\;\hfil}$\crcr
     \noalign{\kern 0pt\nointerlineskip}
     \rightarrowfill\crcr}}\;}}

\def\textlmap#1{\mathop{\vbox{\ialign{
                                 ##\crcr
     ${\scriptstyle\hfil\;\;#1\;\;\hfil}$\crcr
     \noalign{\kern-1pt\nointerlineskip}
     \leftarrowfill\crcr}}\;}}

\newfam\meuffam

\font\tenmeuf=eufm10
\font\sevenmeuf=eufm7
\font\fivemeuf=eufm5

\textfont\meuffam=\tenmeuf
\scriptfont\meuffam=\sevenmeuf
\scriptscriptfont\meuffam=\fivemeuf

\def\germ{\fam\meuffam\tenmeuf}

\def\fg{{\germ f}}
\def\g{{\germ g}}

\def\kg{{\germ k}}

\def\Ag{{\germ A}}

\def\Fg{{\germ F}}

\def\Pg{{\germ P}}

\newtheorem{sz}{Satz}[section]
\newtheorem{thry}[sz]{Theorem}
\newtheorem{pr}[sz]{Proposition}
\newtheorem{re}[sz]{Remark}
\newtheorem{co}[sz]{Corollary}

\newtheorem{lm}[sz]{Lemma}

\begin{document}
\def\Pr{{\rm Pr}}
\def\tr{{\rm Tr}}
\def\End{{\rm End}}
\def\Aut{{\rm Aut}}
\def\Spin{{\rm Spin}}
\def\U{{\rm U}}
\def\SU{{\rm SU}}
\def\SO{{\rm SO}}
\def\PU{{\rm PU}}
\def\GL{{\rm GL}}
\def\spin{{\rm spin}}
\def\su{{\rm su}}
\def\so{{\rm so}}
\def\ub{\underbar}
\def\pu{{\rm pu}}
\def\Pic{{\rm Pic}}
\def\Iso{{\rm Iso}}
\def\NS{{\rm NS}}
\def\deg{{\rm deg}}
\def\Hom{{\rm Hom}}
\def\Aut{{\rm Aut}}
\def\h{{\germ h}}
\def\Herm{{\rm Herm}}
\def\Vol{{\rm Vol}}
\def\pf{{\bf Proof: }}
\def\id{{\rm id}}
\def\i{{\germ i}}
\def\im{{\rm im}}
\def\rk{{\rm rk}}
\def\ad{{\rm ad}}
\def\h{{\bf H}}
\def\coker{{\rm coker}}
\def\dbar{\bar{\partial}}
\def\Lo{{\Lambda_g}}
\def\niq{=\kern-.18cm /\kern.08cm}
\def\Ad{{\rm Ad}}
\def\RSU{\R SU}
\def\ad{{\rm ad}}
\def\dva{\bar\partial_A}
\def\da{\partial_A}
\def\p{\partial\bar\partial}
\def\sp{\Sigma^{+}}
\def\sm{\Sigma^{-}}
\def\spm{\Sigma^{\pm}}
\def\smp{\Sigma^{\mp}}
\def\Tors{{\rm Tors}}
\def\st{{\rm st}}
\def\s{{\rm s}}
\def\oo{{\scriptstyle{\cal O}}}
\def\ooo{{\scriptscriptstyle{\cal O}}}
\def\sw{Seiberg-Witten }
\def\pa{\partial_A\bar\partial_A}
\def\Dr{{\raisebox{0.15ex}{$\not$}}{\hskip -1pt {D}}}
\def\gr{{\scriptscriptstyle|}\hskip -4pt{\g}}
\def\subsetint{{\  {\subset}\hskip -2.45mm{\raisebox{.28ex}
{$\scriptscriptstyle\subset$}}\ }}
\def\ra{\rightarrow}
\def\kod{{\rm kod}}

\title{Donaldson theory on non-K\"ahlerian surfaces and  class $VII$ surfaces with
$b_2=1$}
\author{Andrei Teleman}
\date{May 1, 2005}
\maketitle

\tableofcontents

\begin{abstract}
We prove that any class $VII$ surface with $b_2=1$ has curves. This implies the ``Global
Spherical Shell conjecture"  in the case $b_2=1$: \\

{\it Any minimal class $VII$ surface with $b_2=1$ admits a global
spherical shell, hence it is isomorphic to one of the surfaces in the known list.} 
\\

By the results in \cite{LYZ}, \cite{Te1},  which treat the case $b_2=0$
and give complete proofs  of    Bogomolov's theorem,  one has a complete
classification of all class
$VII$-surfaces with
$b_2\in\{0,1\}$.

 The main idea of  the proof is to show that
a certain moduli space of  $PU(2)$-instantons  on a surface $X$ with no curves (if such a surface existed)
would contain  a closed Riemann surface $Y$ whose  general points correspond to non-filtrable holomorphic
bundles on $X$.  Then we  pass  from a family of bundles on $X$ parameterized by $Y$
to a family of bundles on $Y$ parameterized by $X$, and we use the algebraicity of $Y$
to obtain a contradiction.

The proof   uses essentially techniques from Donaldson theory: compactness theorems for moduli spaces
of
$PU(2)$-instantons   and the Kobaya\-shi-Hitchin correspondence on   surfaces.  
\end{abstract}


\section{Introduction}

\subsection{The main result}

A class $VII$ surface is a compact complex surface $X$ with $b_1(X)=1$ and $\kod(X)=-\infty$. 
The topological invariants of such a surface are 
$$c_2(X)=-c_1(X)^2=b_2(X)\ ,\ b_+^2(X)=0\ .
$$

Therefore the intersection form $H^2(X,\Z)/{\rm Tors}\times H^2(X,\Z)/{\rm Tors}\to\Z$  
is negative definite so, by Donaldson's first theorem, it is trivial over $\Z$.

Class $VII$ surfaces are not classified yet. This is probably the
most important gap in the Enriques-Kodaira classification table.  The case
$b_2=0$ is completely understood: 
\begin{thry}\label{bogo}
Every class $VII$-surface with $b_2=0$ is isomorphic to either a Hopf surface or an Inoue surface.
\end{thry}

This result was stated by Bogomolov a long time ago   \cite{Bo1},
\cite{Bo2}, but his proof is long and difficult to follow (see \cite{BHPV}
p. 230);   complete proofs appeared in
\cite{Te1} and
\cite{LYZ}.

Both proofs are based on the Kobayashi-Hitchin correspondence on non-K\"ahle\-rian surfaces
(see \cite{Bu1}, \cite{LY}, \cite{LT1}) applied to a single holomorphic
bundle: the tangent bundle of the surface.  \\ 

The main result of this paper is:
\begin{thry}\label{main} Let $X$ be a class $VII$ surface with $b_2(X)=1$. Then $X$ has an
effective divisor $C>0$ with
$$c_1^\Q({\cal O}(C)) \in\{\pm c_1^\Q({\cal K}_X),0,2c_1^\Q({\cal K}_X)\}\ ,
$$
where $c_1^\Q$ stands for the first Chern class in rational cohomology.
\end{thry} 

Using Theorem 11.2 in  \cite{Na},  one concludes that
\begin{co} Any minimal class $VII$ surface
$X$ with $b_2(X)=1$ possesses a spherical shell, hence it belongs to the known class of
surfaces.
\end{co}

The proof of the main theorem is again based  on the Kobayashi-Hitchin
correspondence but, unlike   Theorem \ref{bogo} --  which  uses  the 
correspondence for a single bundle --  it requires a careful examination of the
geometry of a certain moduli space of instantons (stable bundles) on $X$; therefore it
is much closer in spirit to  techniques used in Donaldson theory. In particular, one
needs essentially the  Kobayashi-Hitchin correspondence {\it as an isomorphism of
moduli spaces}, and the compactness theorems for moduli spaces of
$PU(2)$ instantons (which, in the non-K\"ahlerian framework, {\it cannot} be obtained
by  complex geometric methods). 

We mention that many arguments can be partially generalized  for class $VII_0$ surfaces
with higher $b_2$ (see \cite{Te2}). By a result of
Dloussky-Oeljeklaus-Toma \cite{DOT}, any class
$VII_0$ surface which has $b_2$ rational curves, contains a
global spherical shell, hence it belongs to the known class. Therefore, 
the classification problem for class
$VII$ surfaces reduces to the  existence of ``sufficiently many" curves.

\subsection{Donaldson Theory revisited} \label{donaldson}

We explain now the gauge theoretical tools 
needed  in the proof (\cite{Do}, \cite{DK}, \cite{LT1}): \\

Let $(M, g)$ be a compact oriented Riemannian 4-manifold and  $E$ be a
Hermitian 2-bundle on $M$. For a  connection $A$ on $E$,  
denote as usually   by   $F_A$ its curvature and by $F_A^0$ its 
trace-free part.

Put
$L:=\det(E)$.  We fix a Hermitian connection $a\in{\cal A}(L)$ and denote by ${\cal
A}_a(E)$ the space of Hermitian connections on $E$ which induce
$a$ on $L$. Such connections are called sometimes ``oriented connections". ${\cal
A}_a(E)$ is  an affine
space over the vector space  
$A^1(su(E))$; the gauge group $SU(E)$ of unitary automorphisms of determinant 1
acts on this affine space naturally.  
The moduli space of  {\it projectively} ASD connections in ${\cal A}_a(E)$ is
$$
{\cal M}_a^{\rm ASD}(E):=\qmod{\{A\in {\cal A}_a(E)|\
(F_A^0)^+=0\}}{SU(E)}\subset\qmod{{\cal A}_a(E)}{SU(E)}=:{\cal B}_a(E)   
$$

Let $P$ be the principal unitary frame bundle of $E$ and  $\bar P:=P/{S^1}$  the
associated $PU(2)$-bundle.  One has a natural identification ${\cal A}_a(E)\simeq{\cal
A}(\bar P)$ which yields a surjection
\begin{equation}
{\cal M}_a^{\rm ASD}(E)\map {\cal M}^{\rm ASD}(\bar P)
\end{equation} 
with finite fibers. This surjection is an isomorphism in the simply connected case, but
in general is {\it not}! The point is that the gauge group $SU(E)$ could be slightly
smaller than the automorphism group ${\rm Aut}(\bar P)$.  This phenomenon can be
easily understood as follows (see \cite{LT1}, p. 141-145 for details):

An  element $\rho\in H^1(X,\Z_2)$ can
be interpreted as a flat $S^1$-connection $a_\rho$ on a Hermitian line bundle
$L_\rho$, which comes with a tautological unitary isomorphism $L_\rho^{\otimes
2}=X\times\C$. The  bundles $E$ and 
$E\otimes L_\rho$  have the same determinant line bundle and are
isomorphic; more precisely, there exists a unitary isomorphism $f:E\to
E\otimes L_\rho$ with $\det(f)\equiv 1$. The map
  $ A \mapsto  f^{-1}(A\otimes a_\rho) $ descends to a well defined map
$\otimes\rho:{\cal B}_a(E)\to{\cal B}_a(E)$. In this way one obtains an 
$H^1(X,\Z_2)$-action  
 on the quotients ${\cal B}_a(E)$, ${\cal M}^{\rm ASD}_a(E)$, and
${\cal M}^{\rm ASD}(\bar P)$ is just the $H^1(X,\Z_2)$-quotient of ${\cal
M}^{\rm ASD}_a(E)$.

Note that the irreducible part $[{\cal
M}_a^{\rm ASD}(E)]^{\rm irr}$ of  ${\cal M}_a^{\rm ASD}(E)$ can contain fixed  points of
this action; in other words there exist in general irreducible  
unitary connections   which project on reducible $PU(2)$-connections.

In classical
gauge theory, one usually works with the simpler moduli space  ${\cal
M}^{\rm ASD}(\bar P)$ and  ignores  ${\cal M}_a^{\rm ASD}(E)$, because   
both spaces should carry equivalent differential topological information.
\vspace{3mm} 
 
{\it For our proof it is    important to consider ${\cal M}_a^{\rm
ASD}(E)$ (rather than ${\cal
M}^{\rm ASD}(\bar P)$), because the 
$H^1(X,\Z_2)$-symmetry of this space plays a crucial role in the proof.}
\\

The precise form of the Kobayashi-Hitchin correspondence we   need is the following
(see \cite{Bu1}, \cite{LT1}, \cite{LT2}, \cite{LY}):
\begin{thry}\label{KH}  Let $(X,g)$ be a compact complex surface endowed with a
Gauduchon metric \cite{G} and
$E$ a Hermitian bundle on $X$. Fix a holomorphic structure ${\cal L}$ on the Hermitian
line bundle $L:=\det(E)$ and let $a\in{\cal A}(L)$ be the corresponding Chern connection.
There is a natural real analytic isomorphism 
$$[{\cal M}^{\rm ASD}_a(E)]^{\rm irr}\textmap{ {KH}\simeq}{\cal M}^{\rm st}_g(E,{\cal
L})\ .$$
\end{thry}

Here ${\cal M}^{\rm st}_g(E,{\cal L})$ denotes the moduli space of $g$-stable
holomorphic structures on $E$ which induce ${\cal L}$ on $L=\det E$, modulo the
complex gauge group ${\rm SL}(E)$  (see \cite{LT1},   \cite{LT2}).

When $h$ is chosen such that $\det(h)$  is a Hermitian-Einstein metric on ${\cal L}$,
this statement follows formally from the standard Kobayashi-Hitchin correspondence
between irreducible Hermitian-Einstein connections and stable bundles 
\cite{LT1}, \cite{LY}. The general statement is a very special case of  the universal
Kobayashi-Hitchin correspondence for oriented pairs    \cite{LT2}, but it
can be easily deduced from the case when   $\det(h)$  is  
Hermitian-Einstein, by noting that the left hand moduli space is in fact
independent of $a$, up to canonical isomorphism. 

Denote by $\lambda$ the semiconnection on $L$ given by the Dolbeault operator of
${\cal L}$. Let ${\cal A}^{0,1}(E)$ be the complex affine space of 
semiconnections (``(0,1)-connections") on
$E$ (see \cite{Do}, \cite{LO}, \cite{LT1}) and ${\cal A}^{0,1}_{\lambda}(E)$  the
subspace of semiconnections which induce
$\lambda$ on $L$.

The Kobayashi-Hitchin isomorphism $KH$ is induced by the affine map 
$${\cal A}_a\to 
{\cal A}^{0,1}_{\cal L}(E)\ ,\ A\mapsto \bar\partial_A\ .$$

Using a standard corollary to Uhlenbeck's compactness theorem \cite{DK},
one gets the following important result, which cannot be obtained by
complex geometric methods.

\begin{co}\label{compactness} In the conditions of Theorem \ref{KH} suppose that
$\rk(E)=2$, and $4c_2(E)-c_1(E)^2\leq 3$. Then the complex moduli space ${\cal
M}^{\rm st}(E,{\cal L})$ can be identified with an open set of the \ub{compact} moduli
space
${\cal M}^{\rm ASD}_a(E)$. Therefore, in this case, ${\cal M}^{\rm
st}(E,{\cal L})$ can be compactified by adding \ub{only} the   reducible
part of ${\cal M}^{\rm ASD}_a(E)$, which can be  identified with the set
of split polystable bundles  ${\cal M}\oplus[{\cal L}\otimes{\cal
M}^{-1}]$, $2\deg_g({\cal M})=\deg_g({\cal L})$, $c_1({\cal
M})(c_1({\cal L})-c_1({\cal M}))=c_2(E)$.   
\end{co}
The compactness of ${\cal M}^{\rm ASD}_a(E)$ follows from the fact that every
$PU(2)$-instanton has  non-positive Pontrjagin number  so, under the assumption
$p_1(\bar P) \geq -3$, the lower strata of the Uhlenbeck compactification of ${\cal
M}^{\rm ASD}(\bar P)$ are empty.
\\ \\
{\bf Remark:} Suppose that we are in the conditions of the above corollary. When
$b_1(X)$ is odd, the compactification
${\cal M}^{\rm ASD}_a(E)$ of
${\cal M}^{\rm st}(E,{\cal L})$ is {\it not} in general a complex space.  
\\  

For instance one can get a moduli space of stable bundles isomorphic to an open disk
$D\subset\C$, which is compactified in the natural way by adding a circle.
The point is that the stratum of  reducible connections in ${\cal M}^{\rm ASD}_a(E)$
(split polystable bundles  with fixed determinant) can have odd real dimension.

This remark also shows that, in the non-K\"ahlerian framework,  one cannot hope to find
a purely complex geometric way to compactify the moduli spaces of stable bundles (as
the Gieseker compactification  in the algebraic
case).
\\

The Kobayashi-Hitchin correspondence (Theorem \ref{KH}) can be regarded as a practical method for
computing moduli spaces of projectively ASD connections with complex geometric methods. Note however that the classification of holomorphic bundles on a non-algebraic manifold  is in general a very difficult
problem. Indeed, on a non-algebraic manifold there exists in general   non-filtrable   bundles (see
section \ref{clsimple}). Such bundles are always stable  with respect to any Gauduchon
metric. On the other hand there exists no general construction or classification method for
non-filtrable bundles \cite{BLP}.  
\subsection{The strategy of the proof}\label{strategy}

The idea of the proof is the following: \\

Let $X$ be a class $VII$ surface with $b_2(X)=1$. We fix a Gauduchon metric $g$ on $X$
such that $\deg_g({\cal K}_X)\ne 0$ (which is possible since $c_1({\cal K}_X)^2\ne 0$), and
consider the moduli space
${\cal M}^{\rm st}:= {\cal M}^{\rm st}_g(0,{\cal K}_X)$ of $g$-stable bundles ${\cal E}$
with
$\det({\cal E})\simeq{\cal K}_X$ and $c_2=0$. The expected complex dimension of this space
is 1. 

Put
$$p:=|\Tors(H^2(X,\Z))|\ ,\ q:=|\Tors_2(H^2(X,\Z))|\ ,
$$
where $\Tors_2$ stands for the 2-torsion of an abelian group.
One has $0<q\leq p$.   Supposing that $X$ has no curves $C$ with 
$$c_1^\Q({\cal O}(C)) \in\{\pm c_1^\Q({\cal K}_X),0,2c_1^\Q({\cal K}_X)\} \ ,
$$
 we show that  the compactification $\overline{{\cal M}^{\rm st}}$  of
${\cal M}^{\rm st}$ given by Corollary \ref{compactness} contains {\it as an open set}
the   disjoint union
$$\left[\coprod_{i=1}^{p-q} D_i\right]\coprod \left[\coprod_{j=1}^{q} P_j\right]
$$
where $D_i$ are closed disks, and $P_j$ are copies  of     one of the  two spaces
 illustrated below
$$
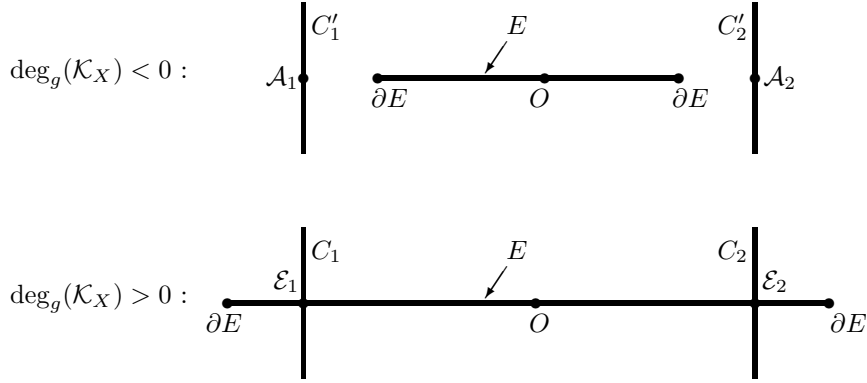
\begin{array}{c}
\unitlength=1mm
\begin{picture}(90,25)(-55,-15)
\linethickness{0,6mm}  
\put(-69,0){$\deg_g({\cal K}_X)<0$ :}
\put(-20,0){\line(1,0){40}} 
 \put(0,-0,9){ $\bullet$ }
 \put(0,-3,9){$O$}
 \put(-3,5,9){$E$}
 \put(-3,4,9){\vector(-2,-3){2.8}}
 \put(-29,5,9){$C'_1$} 
\put(25,5,9){$C'_2$} 
 \put(-21,-0,9){$\bullet$}
 \put(-21,-3,9){$\partial E$}
\put(19,-0,9){$\bullet$}
\put(19,-3,9){$\partial E$}
\put(-30,-10){\line(0,1){20}}
\put( 30,-10){\line(0,1){20}}
 \put(-30.9,-0,9){$\bullet$}
 \put(-35,-0,9){${\cal A}_1$} 
 \put(29.15,-0,9){$\bullet$} 
 \put(31.15,-0,9){${\cal A}_2$}
 \end{picture} 
\end{array}   
$$
$$
\begin{array}{c}
\unitlength=1mm
\begin{picture}(90,25)(-55,-15)
\linethickness{0,6mm}  
\put(-69,0){$\deg_g({\cal K}_X)>0$ :}
\put(-40,0){\line(1,0){80}}
 \put(0,-0,9){$\bullet$}
 \put(-41,-0,9){$\bullet$}
 \put(-34,1,9){${\cal E}_1$}
\put(31,1,9){${\cal E}_2$}
 \put(-29,5,9){$C_1$} 
\put(25,5,9){$C_2$}
\put(-3,5,9){$E$}
 \put(-3,4,9){\vector(-2,-3){2.8}}
\put(39,-0,9){$\bullet$}
\put(-30,-10){\line(0,1){20}}
\put( 30,-10){\line(0,1){20}}
 \put(-30.9,-0,9){$\bullet$}
 \put(29.15,-0,9){$\bullet$}
 \put(-43,-3,9){$\partial E$}
 \put(40,-3,9){$\partial E$}
\put(0,-3,9){$O$}
 \end{picture} 
\end{array}   
$$
 One must choose the first picture when
$\deg_g({\cal K}_X)<0$ and the second when
$\deg_g({\cal K}_X)>0$. In both pictures, the segment $E$ represents a closed disk with
center $O$. The points of $\cringle{E}\setminus\{O\}$ correspond to   filtrable  stable
bundles, which can be completely classified. The   origin $O$ is a non-filtrable bundle
which is obtained as pushforward of  a  holomorphic line bundle on a
bicovering; it is a fixed   point under a natural $\Z_2$-action. The circle
$\partial E$ is the stratum of reducible solutions  (split polystable  bundles)
mentioned in Corollary
\ref{compactness}. Moreover, the locus of reducible  solutions  is the union
of the  boundaries $\partial E_j\subset P_j$  and   the boundaries  
  $\partial D_i$ (hence it is a disjoint union of $p$ circles). 

Each vertical segment  $C'_i$  (respectively $C_i$) represents a small open disk 
whose only filtrable point is the center ${\cal A}_i$ (respectively ${\cal
E}_i$). The   ${\cal E}_i$'s are the only singularities of the moduli space in
the second case, whereas   ${\cal M}^{\st}$ is smooth in the first case. 

The main idea is the following: since $\overline{{\cal M}^{\rm st}}$ is compact, 
it follows that the disks $C'_i$ ($C_i$) are part of a  compact complex subspace of
the moduli space. One can easily see that this compact subspace is in fact a smooth
(possibly non-connected) closed Riemann surface.

The crucial point here is  that our disjoint union is embedded {\it injectively}
as an open set in the moduli space, so that the vertical and the horizontal loci belong to
different irreducible components.

The presence in the moduli space of a smooth closed complex curve $Y$  which
contains both filtrable and  non-filtrable points   leads  to a
contradiction (see section \ref{families}). The first step is to   put together the bundles $({\cal
E}_y)_{y\in Y}$ and get a bundle
${\cal E}$ on
$Y\times X$. Since $H^2(Y,{\cal O}_Y^*)=0$, there is no obstruction to the construction
of such a classifying family ${\cal E}$.  In particular, one has a family $({\cal E}^x)_{x\in X}$ of
bundles on the curve
$Y$ parameterized by $X$. We explain briefly, in two simple particular
situations, why such a family cannot exist: First, when  the 
bundles
${\cal E}^x$ are all  semistable, one would get a morphism $\varepsilon$
from $X$ into a  moduli space of  semistable bundles over $Y$, {\it which is an
algebraic geometric object}.  Since our surface has algebraic dimension 0, it follows
that only the case
$\varepsilon$=constant (which gives a very simple family) must be eliminated.

The   opposite case is when  the  $({\cal E}^x)_{x\in X}$ are
non-semistable for all $x\in X$. A non semistable 2-bundle ${\cal F}$
over a curve has a unique maximal destabilizing sub-line bundle ${\cal D}({\cal
F})$, which depends ``meromorphically" on ${\cal F}$. Since the components of
$\Pic(Y)$ are also projective manifolds, one has only to discuss the case when
the map
$x\mapsto [{\cal D}({\cal E}^x)]$ is constant. Putting together the maximal
destabilizing line bundles ${\cal D}({\cal E}^x)$ we get a subsheaf of
rank 1 of ${\cal E}$, contradicting the fact  that ${\cal E}_y$ is
non-filtrable for generic $y\in Y$.

\section{Background material}

\subsection{The Picard group,  the Gauduchon degree and the square roots
of $[{\cal O}]$}\label{roots}

Let $X$ be a surface, and $\Pic(X)=H^1(X,{\cal O}^*)$ its Picard group.    For a
cohomology class $c\in NS(X)$, we will denote by $\Pic^c(X)$ the  corresponding
connected component of $\Pic(X)$. $\Pic^0(X)$ is  just the identity component of the group
$\Pic(X)$. Let 
$\Pic^T(X)$ be the subgroup of line bundles with torsion Chern class. Therefore
$$\qmod{\Pic(X)}{\Pic^0(X)}\simeq NS(X) \ ,\ \qmod{\Pic(X)}{\Pic^T(X)}\simeq \qmod{NS(X)}{{\rm
Tors}(NS(X))}\ .
$$

Let $g$ be a Gauduchon metric \cite{G} on  $X$, i.e. a Hermitian metric such that
$\partial\bar\partial\omega_g=0$. The degree map   associated with
$g$ is  a   group morphism
$$\deg_g:\Pic(X)\map\R
$$
defined by
$$\deg([{\cal L}]):=\int_X c_1({\cal L},h)\wedge \omega_g\ ,
$$
where $c_1({\cal L},h)$ denotes the Chern {\it form} of any Hermitian metric on ${\cal L}$. The degree
map is a topological invariant (i.e. it vanishes on $\Pic^T(X)$) if and only if 
$b_1(X)$ is even (see
\cite{LT1}).   The restriction
$\resto{\deg_g}{\Pic^T(X)}$ is independent of $g$ up to {\it positive} multiplicative constant
(\cite{LT1}, p. 41).  Its kernel
$$\Pic^{f}(X):=\ker(\resto{\deg_g}{\Pic^T(X)}:\Pic^T(X)\to\R)\subset \Pic^T(X)
$$
is always compact, because the Kobayashi-Hitchin correspondence gives an isomorphism
of real Lie groups
\begin{equation}\label{flat}
\Hom(\pi_1(X),S^1)\simeq\Pic^{f}(X) 
\end{equation}
 (see \cite{LT1}). In the case of
surfaces with even
$b_1$  one has
$\Pic^{  f}(X)=\Pic^T(X)$, whereas  for surfaces with odd $b_1$, $\Pic^{  f}(X)$ is a {\it
real} codimension 1 compact subgroup of $\Pic^T(X)$. 

Let $X$ be an arbitrary  class $VII$ surface. For such a
surface the identity-connected component $\Pic^0(X)$ of  $\Pic(X)$ is isomorphic to
$\C^*$.  More precisely,   the natural morphisms 
$$H^1(X,\C^*)\map \Pic^T(X)\ ,\ \qmod{H^1(X,\C)}{H^1(X,\Z)}\map\Pic^0(X)
$$
are isomorphisms. Therefore, any holomorphic line bundle ${\cal L}$  with torsion Chern
class has a unique flat connection, and this flat connection is Hermitian if and only if
the degree of
${\cal L}$ (with respect to any Gauduchon metric) vanishes.

An important role in this article will be played by {\it the square roots of} the
class $[{\cal O}]$ of the trivial holomorphic line bundle. Using (\ref{flat}), one obtains an
isomorphism 
\begin{equation}\label{iso}
H^1(X,\Z_2)=\Hom(\pi_1(X),\Z_2)\simeq {\rm Tors}_2(\Pic(X)) 
\end{equation}

The long exact sequence associated with $0\to \Z\textmap{2\cdot} \Z\to\Z_2\to 0$ gives
$$0\map \qmod{H^1(X,\Z)}{2H^1(X,\Z)}\map H^1(X,\Z_2)\map {\rm Tors}_2(H^2(X,\Z))\map 0\ ,
$$
where the epimorphism on the right coincides with the Chern class morphism via the isomorphism
(\ref{iso}), and the quotient  on the left can be identified with $\Z_2$ in our case (because $b_1(X)=1)$.

Therefore, for class $VII$ surfaces one has:
\begin{re} A connected component $\Pic^c(X)$ of $\Pic(X)$ contains square roots of 
$[{\cal O}]$ if and only if $c\in {\rm Tors}_2(H^2(X,\Z))$; in this case it contains
two square roots which are conjugate under the natural action of $\Z_2$.
\end{re}

\section{A moduli space of simple bundles} 
 
\subsection{Classifying simple filtrable bundles}\label{clsimple}

Let $X$ be a class $VII$ surface with $b_2=1$. We will use the following
simplified notations:
$${\cal K}:={\cal K}_X\ ,\ \Pic:=\Pic(X)\ ,\ \Pic^c:=\Pic^c(X)\ ,\ \Pic^T:=\Pic^T(X) \
.
$$
Let ${\cal M}^s$ be the moduli space
of simple holomorphic rank 2 bundles ${\cal E}$ with $\det{\cal E}\simeq {\cal K}$
and
$c_2=0$.

We recall that a rank 2 bundle ${\cal E}$ on a complex surface $S$ is called {\it
filtrable} if one of the  following equivalent conditions is satisfied:
\begin{enumerate}
\item ${\cal E}$ has a subsheaf of rank 1.
\item There exists a holomorphic line bundle ${\cal L}$ on $S$ such that $H^0({\cal
L}^\vee\otimes {\cal E})\ne 0$.
\item There exist line bundles ${\cal L}'$ and ${\cal L}''$ on $S$, a dimension
0 locally complete intersection $Z\subset S$ and a short exact sequence of the  form
$$0\map{\cal L}'\map{\cal E}\map {\cal L}''\otimes{\cal I}_Z\map 0\ .
$$
\end{enumerate}

\begin{lm}\label{generator} Let $X$ be a class $VII$ surface with $b_2=1$. Then 
$c_1({\cal K})$ modulo torsion  is a
$\Z$-generator of the rank one $\Z$-module $H^2(X,\Z)/{\rm Tors}$.
\end{lm}

This follows from the well-known formula:
$$c_1({\cal K})^2=c_1(X)^2=-b_2(X)\ ,$$ 
for a class $VII$-surface $X$.\\

\begin{pr}\label{filtrable} Let ${\cal E}$ be a \ub{filtrable} holomorphic rank 2 bundle on
$X$ with
$c_2({\cal E})=0$,
$c_1^\Q({\cal E})=c_1^\Q({\cal K})$. Then there exist  holomorphic   line bundles ${\cal L}$,
${\cal M}$ on
$X$ such that
\begin{enumerate}
\item $c_1^\Q({\cal L})= 0$.
\item $c_1^\Q({\cal M})=c_1^\Q({\cal K})$.
\item ${\cal E}$ is either an extension of ${\cal L}$ by ${\cal M}$ or an extension
of ${\cal M}$ by ${\cal L}$.
\end{enumerate}
\end{pr}
\pf  Choose an exact sequence of the form
$$0\map{\cal L}'\map{\cal E}\map {\cal L}''\otimes{\cal I}_Z\map 0\ .
$$
as above.  By Lemma \ref{generator}, one can write $c_1^\Q({\cal
L'})=nc_1^\Q({\cal K})$, with
$n\in\Z$.
$$c_2({\cal E})=0=|Z|+c_1^\Q({\cal L'}) \cup c_1^\Q({\cal
L''})= |Z|+n(n-1)\ ,
$$
which can only hold when $|Z|=0$,  and   $n\in\{0,1\}$.
\qed
\begin{pr}\label{bundles} Suppose that $X$ has no effective divisor $C>0$
with 
$$c_1^\Q({\cal O}(C)) \in\{\pm c_1^\Q({\cal K}),0,2c_1^\Q({\cal K})\}\ .
$$
\begin{enumerate}
\item For every line bundle ${\cal L}$ with torsion Chern class there exists a
unique (up to isomorphism) rank two  bundle ${\cal E}_{\cal L}$ which is the
central term of a  {nontrivial} extension
\begin{equation}\label{EL}
0\map {\cal L}\map {\cal E}_{\cal L}\map {\cal K}\otimes{\cal
L}^{\vee}\map 0\ .
\end{equation}

For every square root ${\cal R}$ of ${\cal O}$ there exists a
unique (up to isomorphism) rank two   bundle ${\cal A}_{\cal R}$ which is the
central term of a  {nontrivial} extension
\begin{equation}\label{AR}
0\map {\cal R}\otimes{\cal K}\map {\cal A}_{\cal R}\map {\cal R}\map 0\ .
\end{equation}
\item The bundles ${\cal E}_{\cal L}$, ${\cal A}_{\cal R}$ are simple.
Moreover, 
${\cal E}_{\cal L'}\not\simeq  {\cal E}_{\cal L''}$ when
${\cal L}'\not\simeq{\cal L}''$ and ${\cal A}_{\cal R'}\not\simeq  {\cal
A}_{\cal R''}$ when ${\cal R}'\not\simeq{\cal R}''$. 
\item  ${\cal A}_{\cal R}\not\simeq{\cal E}_{\cal L}$, for every $ [{\cal
L}]\in
\Pic^T$ and square root ${\cal R}$ of ${\cal O}$.
\item $[{\cal E}_{\cal L}]$ is a smooth point of ${\cal M}^{\rm s}$,
except  when $[{\cal L}]\in\Tors_2(\Pic)$.
\item For every $[{\cal R}]\in\Tors_2(\Pic)$, ${\cal M}^{\rm s}$ is  reducible
(hence singular) at
$[{\cal E}_{\cal R}]$; in a neighborhood of this point
${\cal M}^{\rm s}$ consists of two smooth curves $C_{\cal R}$ and
$\Phi_{\cal R}$ intersecting transversally at
$[{\cal E}_{\cal R}]$.
$\Phi_{\cal R}$ is just a neighborhood of $[{\cal E}_{\cal R}]$ in   the
1-parameter family $\{{\cal E}_{\cal L}|\ [{\cal L}]\in
\Pic^T\}$.
\item For any    $[{\cal R}]\in\Tors_2(\Pic)$, the moduli space ${\cal M}^{\rm s} $ is smooth at 
$[{\cal A}_{\cal R}]$.  
\item The map 
$$\Fg:\Pic^T\coprod\Tors_2(\Pic)  \map {\cal M}^{\rm s} 
$$
given by  ${\cal L}\mapsto[{\cal E}_{\cal L}]$, ${\cal R}\mapsto[{\cal A}_{\cal
R}]$ parameterizes \ub{bijectively} the filtrable part of
${\cal M}^{\rm s}$. 
\end{enumerate}
\end{pr}
\pf An effective divisor $C\subset X$ with $c_1^\Q({\cal O}(C))= nc_1^\Q({\cal K})$ will be called
numerically $n$-canonical.\\ \\ 
1.  By Riemann-Roch Theorem, one has
$$\chi ({\cal K}^\vee\otimes{\cal L}^{\otimes 2})=-1
$$
and $H^0({\cal K}^\vee\otimes{\cal L}^{\otimes 2})=H^2({\cal K}^\vee\otimes{\cal
L}^{\otimes 2})=0$, if $X$ has no numerically anticanonical (respectively
numerically bicanonical) curves. Therefore
$${\rm Ext}^1({\cal K}\otimes{\cal L}^\vee,{\cal L})=H^1({\cal K}^\vee\otimes{\cal L}^{\otimes 2})
\simeq\C\ .
$$
Therefore, up to  $\C^*$-equivalence, one has a unique nontrivial extension of ${\cal
K}\otimes {{\cal L}}^\vee$ by ${\cal L}$.  For the second type of extensions, note that
$H^1({\cal K})\simeq H^1({\cal O})^\vee\simeq\C$.
\\ \\
2.   A morphism   $\varphi:{\cal
E}_{\cal L'}\to {\cal E}_{\cal L''}$ defines a diagram
$$
\begin{array}{ccccccccc}
0&\map&{\cal L'}&\textmap{\alpha'}& {\cal E}_{\cal L'}&\textmap{\beta'}&{\cal
L'}^{\vee}\otimes{\cal K}&\map& 0
\\
&&&&\downarrow\varphi
\\
0&\map&{\cal L''}&\textmap{\alpha''}& {\cal E}_{\cal L''}
&\textmap{\beta''}&{\cal L''}^{\vee}\otimes{\cal K}&\map& 0
\end{array}
$$
The composition $\beta''\circ\varphi\circ\alpha'$ vanishes because $X$ has no
numerically canonical curves. Therefore, $\beta''\circ\varphi$ induces a morphism $q_\varphi:{\cal
L'}^{\vee}\otimes{\cal K}\to{\cal
L''}^{\vee}\otimes{\cal K}$.
\\ \\
{\it Case a.} ${\cal
L}'\not\simeq{\cal L}''$

In this case $q_\varphi$ vanishes because  $X$ has no numerically trivial
curves. This shows  that $\varphi$ factorizes as $\varphi=\alpha''\circ \psi$, for  a
morphism
$\psi:{\cal E}_{\cal L'}\to {\cal L}''$, so $\varphi$ cannot be an isomorphism. 
Therefore ${\cal E}_{{\cal L}'}\not\simeq {\cal E}_{{\cal L}''}$.
\\
\\
{\it Case b.} ${\cal
L}'={\cal L}''$

In this  case  $q_\varphi$ has the form  $q_\varphi=\zeta\id_{{\cal
L'}^{\vee}\otimes{\cal K}}$, so that
$\varphi_\zeta:=\varphi-\zeta\id_{{\cal E}_{\cal L'}}$  factorizes   as
$\varphi_\zeta=\alpha''\circ \psi_\zeta$, for  a morphism
$\psi_\zeta:{\cal E}_{\cal L'}\to {\cal L}'$. The composition  $\psi_\zeta\circ
\alpha':{\cal L}'\to{\cal L}'$  is trivial, because otherwise $\psi_\zeta$
would define a splitting of the first exact sequence. Therefore
$\psi_\zeta\circ
\alpha'=0$, which shows that $\psi_\zeta$   factorizes through a morphism
${{\cal L}'}^\vee\otimes{\cal K}\to {\cal L}'$. This must vanish,
because  $X$ has no numerical anti-canonical curves.  Therefore
$\varphi=\zeta\id_{{\cal E}_{\cal L'}}$, proving that
${\cal E}_{\cal L'}$ is simple.\\

The same method applies for the statements concerning the bundles ${\cal A}_{\cal R}$. 
\\ \\
3. Let $\varphi:{\cal E}_L\ra {\cal A}_{\cal R}$ be a morphism, and consider the
diagram.
$$
\begin{array}{ccccccccc}
0&\map&{\cal L}&\textmap{\alpha}& {\cal E}_L&\textmap{\beta}&{\cal
L}^{\vee}\otimes{\cal K}&\map& 0
\\
&&&&\downarrow\varphi
\\
0&\map&{\cal K}\otimes{\cal R} &\textmap{a}&{\cal
A}_{\cal R} &\textmap{b}&{\cal R}&\map& 0\ .
\end{array}
$$

The induced morphism $b\circ \varphi\circ \alpha$ is a section of the
holomorphic line bundle ${\cal L}^\vee\otimes{\cal R}$. Since $X$ has no
numerically trivial curves, this morphism can only be trivial if ${\cal
L}={\cal R}$, and in this case,
 it  is a multiple $\zeta\id_{\cal R}$ of the identity map. If $\zeta\ne
0$,
$\varphi\circ\alpha$ will provide a splitting of the  second
line. Therefore $b\circ \varphi\circ \alpha=0$, hence $\varphi$ induces a
morphism
${\cal L}^{\vee}\otimes{\cal K}\to {\cal R}$, which will vanish when
$X$ admits no numerically anticanonical curves.

This shows that $\varphi$ factorizes as $\varphi=a\circ \psi$ for a morphism
$\psi:{\cal E}_L\to {\cal K}\otimes{\cal R}$, so
 it cannot be an isomorphism.  
\\ \\
4. An element $\varphi\in H^0({\cal E}nd_0({\cal E}_{\cal L})\otimes{\cal
K})=H^2({\cal E}nd_0({\cal E}_{\cal L}))^\vee$ defines a morphism $\varphi:{\cal
E}_{\cal L}\to{\cal E}_{\cal L}\otimes{\cal K}$. 

Consider the diagram
$$
\begin{array}{ccccccccc}
0&\map&{\cal L}&\textmap{\alpha}& {\cal E}_{\cal L}&\textmap{\beta}&{\cal
L}^{\vee}\otimes{\cal K}&\map& 0
\\
&&&&\downarrow\varphi
\\
0&\map&{\cal L}\otimes{\cal K}&\textmap{\alpha\otimes\id}& {\cal E}_{\cal
L}\otimes{\cal K}&\textmap{\beta\otimes\id}&{\cal L}^{\vee}\otimes{\cal K}^{\otimes
2}&\map& 0
\end{array}
$$
 When $X$ has no bicanonical divisors, one has
$(\beta\otimes\id)\circ\varphi\circ\alpha=0$, so $\varphi$ induces a morphism 
${\cal
L}^{\vee}\otimes{\cal K}\to {\cal
L}^{\vee}\otimes{\cal K}^{\otimes 2}$, which must vanish, because $H^0({\cal K})=0$.
Therefore, $\varphi$ factorizes as $\varphi=(\alpha\otimes\id)\circ \psi$  for  a
morphism $\psi:{\cal E}_{\cal L}\to {\cal L}\otimes {\cal K}$.  The composition
$\psi\circ \alpha$ must vanish, because $X$ has no canonical divisors.

Therefore $\psi$ is induced by a morphism ${\cal
L}^{\vee}\otimes{\cal K}\to {\cal
L} \otimes{\cal K}$ which can be nontrivial only when ${\cal L}^{\otimes 2}\simeq
{\cal O}$.
\\ \\  
5. We will  study the Kuranishi local model of ${\cal M}^s$ at the point  $[{\cal E}_{\cal R}]$.

The bundle ${\cal E}nd_0({\cal E}_{\cal R})$ fits in the following diagram with exact horizontal and
vertical lines. 
$$\begin{array}{ccccccccc}
&&0&&&\\ 
&&\downarrow&&&\\ 
&&\ {\cal K}^\vee&&&&&&\\ 
&&\ \downarrow s&&&&&&\\ 
0&\map&{\cal U}&\stackrel{j}{\hookrightarrow}&{\cal E}nd_0({\cal E}_{\cal R})&\textmap{h}&{\cal
K}&\map&0\ .\\  &&\ \downarrow t&&&&&&\\ 
&& {\cal O}&&&&&&\\ 
&& \downarrow&&&&&&\\ 
&& 0&&&&&& 
\end{array}
$$
Here ${\cal U}$ is just the kernel of the morphism $h:{\cal E}nd_0({\cal E}_{\cal R})\to{\cal K}$ given by
$$h:\varphi\mapsto \beta\circ\varphi\circ\alpha$$
 (with the notations of 1. and 2.), i.e. the sheaf of trace free
endomorphisms of ${\cal E}_{\cal R}$  which leave invariant the
filtration  $0\subset  \alpha({\cal R})\subset{\cal E}_{\cal R}$. The
fiber
${\cal U}(x)$ is just the Lie algebra of the parabolic subgroup of $SL({\cal E}_{\cal R}(x))$ which is the
stabilizer of the line $\alpha({\cal R})(x)\subset{\cal E}_{\cal R}(x)$.

The morphism 
 $s:{\cal
K}^\vee={\cal H}om({\cal K}\otimes{\cal R},{\cal R})\map {\cal U}$ 
 is given by 
$$\psi\mapsto \alpha\circ\psi\circ\beta\ ,$$
 and $t$ is defined by 
$$u|_{\cal R}=t(u)\id_{\cal
R}\ .$$

 The vertical extension is non-trivial (it is just the extension defining ${\cal
E}_{\cal R}$ tensorized by ${\cal K}^\vee\otimes{\cal R}$), so one gets immediately 
$$H^1({\cal U})=H^1({\cal
O})\simeq\C\ ,\
 H^2({\cal U})=0\ .$$  
The long exact sequence associated with the horizontal line gives an exact sequence
\begin{equation}\label{seq}
0\map H^1({\cal U})=H^1({\cal O})\map H^1({\cal E}nd_0({\cal E}_{\cal R}))\map
H^1({\cal K})\map 0
\end{equation}
and $H^2({\cal E}nd_0({\cal E}_{\cal R}))\simeq H^2({\cal K})\simeq\C$.

The germ of ${\cal M}^s$ at $[{\cal E}_{\cal R}]$ is isomorphic to the vanishing locus of a holomorphic
map $H^1({\cal E}nd_0({\cal E}_{\cal R}))\supset V\textmap{\chi} H^2({\cal E}nd_0({\cal E}_{\cal R}))$,
defined in a neighborhood of $V$ of 0, and having the properties 
$$\chi(0)=0,\ d_0\chi=0\ ,\ d^{(2)}_0\chi(u,v)=[u, v]
$$
where $[\cdot,\cdot]$ is the symmetric bilinear map 
$$H^1({\cal E}nd_0({\cal E}_{\cal R}))\times H^1({\cal
E}nd_0({\cal E}_{\cal R}))\to H^2({\cal E}nd_0({\cal E}_{\cal R}))$$
 induced by the commutator map on trace free
endomorphisms.

Using a fiber splitting ${\cal E}_{\cal R}(x)={\cal R}(x)\oplus[{\cal K}\otimes{\cal R}](x)$	of the exact
sequence (\ref{EL}), one has
$$\left[\left(\matrix{\zeta&\lambda\cr k&-\zeta}\right),
\left(\matrix{z&l\cr0&-z}\right)\right]=\left(\matrix{*&*\cr 2 k z&*}\right)\ ,\
\forall z,\ \zeta\in\C,\ l,\ \lambda\in {\cal K}^{\vee}(x)\ ,\ k\in{\cal K}(x)\ .
$$
This shows that, in any point
$x\in X$, it holds
$$h([\varphi,u])= 2 h(\varphi)t(u)\ ,  \forall\varphi\in  {\cal
E}nd_0({\cal E}_{\cal R})_x\ u\in {\cal U}_x.
$$
Therefore, via the isomorphism  $h_*:H^2({\cal
E}nd_0({\cal E}_{\cal R}))\to H^2({\cal K})$,   one can write 
$$[\varphi,j_*(u)]=2h_*(\varphi)\cup t_*(u) \ ,  \forall\varphi\in H^1({\cal E}nd_0({\cal E}_{\cal R})),\
\forall u\in H^1({\cal U})\  . $$
Via an isomorphism $H^1({\cal E}nd_0({\cal E}_{\cal
R}))=H^1({\cal K})\oplus H^1({\cal O})\simeq\C^2$ defined by a splitting of (\ref{seq}), the quadratic
form associated with the second derivative $d^{(2)}_0\chi$ will have the form
$$(k,\zeta)\mapsto 2k \zeta+ \alpha k^2\ .
$$
Since the first derivative vanishes, it follows that the vanishing locus $Z(\chi)$ has a simple normal
crossing singularity at $0$.

It remains to prove that the map ${\cal L}\mapsto{\cal E}_{\cal L}$ parameterizes a  curve in the
moduli space whose tangent space at $[{\cal E}_{\cal R}]$ is $j_*(H^1({\cal U}))\subset H^1({\cal
E}nd_0({\cal E}_{\cal R}))$.

Let $\Theta$ be the underlying  differentiable line bundle of ${\cal R}$, let $\theta:=c_1(\Theta)\in{\rm
Tors}_2(H^2(X,\Z))$, and $[{\cal L}]\in\Pic^\theta$. Let $K$ be the underlying
${\cal C}^\infty$ differentiable bundle of ${\cal K}$. The bundle
${\cal E}_{\cal L}$ is obtained by putting on the differentiable bundle
$\Theta\oplus [K\otimes\Theta]$ the holomorphic structure defined by an integrable
semiconnection of the form
$$\delta=\left(\matrix{ \delta'&\sigma\cr 0&\delta''}\right)\in{\cal A}^{0,1}(\Theta\oplus
[K\otimes\Theta])\ .
$$
where $\delta'$ defines the holomorphic structure ${\cal L}$ on  
$\Theta$, $ \delta''$ defines the holomorphic structure ${\cal K}\otimes{\cal L}^\vee$ on
$K\otimes\Theta\simeq K\otimes\Theta^{-1}$ and
$\sigma\in A^{0,1}(K^\vee)$ satisfies 
\begin{itemize}
\item[i)] $\delta'\circ \sigma +\sigma\circ \delta''=0$, i.e. $\sigma$ is
$\bar\partial$-closed with respect to the holomorphic structure
$\delta'\otimes(\delta'')^\vee$ (which defines the holomorphic structure
${\cal K}^\vee\otimes {\cal L}^{\otimes 2}$  on $K^\vee$).
\item[ii)] The Dolbeault $\bar\partial$-cohomology class defined by $\sigma$ is non-zero. 
\end{itemize}

Let $S$ be a small neighborhood of  $l_0:=[{\cal R}]$ in $\Pic^\theta$
and $S\ni l\mapsto\delta'_l$ a holomorphic family of integrable
semiconnections (``(0,1)-connections")  on
$\Theta$, such that, for any $l$,  the isomorphism class of the
holomorphic line bundle 
${\cal L}_l:=(\Theta,\delta'_l)$ is $l$. Put $\delta''_l:=\kappa\otimes(\delta'_l)^\vee$,
where $\kappa$ is the semiconnection on $K$ which defines the canonical holomorphic structure
${\cal K}$. The point is that, for sufficiently small $S$, we can choose  $\sigma_l$
satisfying properties i), ii) above for
$(\delta'_l,\delta''_l)$ and depending holomorphically on $l\in V$.  Indeed, since $h^1({\cal
K}^\vee\otimes {\cal L}_l^{\otimes 2})$ is constant, the family of kernels of the operators
 $\delta'_l\otimes(\delta'')^\vee_l: A^{0,1}(K^\vee)\to A^{0,2}(K^\vee)$ 
gives a holomorphic map from
$S$ to the Grassmann manifold of closed subspaces of (a suitable Sobolev
completion of) $A^{0,1}(K^\vee)$.

Let $v\in T_{l_0}(\Pic^\theta)\setminus\{0\}$, and consider the map 
$\delta:S\to{\cal A}^{0,1}_\kappa(\Theta\oplus [K\otimes\Theta])$

$$l\mapsto \delta_l:=\left(\matrix{ \delta'_l&\sigma_l\cr 0&\delta''_l}\right)\ .
$$
The derivative $d_{l_0}\delta(v)$ is  obviously an element in 
$A^{0,1}({\cal U})$   which is
$\bar\partial$-closed   with respect to the holomorphic structure $\delta_{l_0}$
(i.e. the holomorphic structure induced from
${\cal E}nd_0({\cal
E}_{\cal R})$), so it defines a Dolbeault cohomology class $[d_{l_0}\delta(v)]\in H^1({\cal U})$.
Therefore the corresponding tangent vector
$w\in T_{[{\cal E}_{\cal R}]}({\cal M}^s)=H^1({\cal E}nd_0({\cal
E}_{\cal R}))$ belongs indeed to  $j_*(H^1({\cal U}))$. Moreover, this 
vector does not vanish    because, differentiating   the identity
$[(\Theta,\delta'_l)]=l$  gives 
$t_*([d_{l_0}\delta(v)])=v$.  
\\ \\
 6.   An element $\varphi\in H^0({\cal E}nd_0({\cal A}_{\cal R})\otimes{\cal
K})=H^2({\cal E}nd_0({\cal A}_{\cal R}))^\vee$ defines a morphism $\varphi:{\cal
A}_{\cal R}\to{\cal A}_{\cal R}\otimes{\cal K}$. 

Consider the diagram
$$
\begin{array}{ccccccccc}
0&\map&{\cal R}\otimes {\cal K}&\textmap{\alpha}& {\cal A}_{\cal R}&\textmap{\beta}&{\cal R}&\map& 0
\\
&&&&\downarrow\varphi
\\
0&\map&{\cal R}\otimes{\cal K}^{\otimes 2}&\textmap{\alpha\otimes\id}& {\cal A}_{\cal
R}\otimes{\cal K}&\textmap{\beta\otimes\id}&{\cal R}\otimes{\cal K} &\map& 0
\end{array}
$$
One has
$(\beta\otimes\id)\circ\varphi\circ\alpha=0$, because, otherwise, one would get a  
splitting of the second exact sequence. One proceeds as in {\it Case b.} above, taking
into account that
$X$ has no canonical or bicanonical curves, and get $\varphi=0$.
\\ \\
7. First of all note that an extension of the type
$$0\map {\cal M}\map {\cal E}\map {\cal K}\otimes{\cal M}^{-1}\map 0
$$
with $c_1^\Q({\cal M})=c_1^\Q({\cal K})$ can be nontrivial only if the
line bundle
${\cal M}$ has the form ${\cal M}={\cal K}\otimes{\cal R}$ with $[{\cal
R}]\in\Tors_2(\Pic)$. Indeed, $\chi({\cal K}^\vee\otimes{\cal
M}^{\otimes2})=0$ by Riemann Roch, and $h^0({\cal K}^\vee\otimes{\cal
M}^{\otimes2})=0$ because $X$ has no numerically canonical curves.
Therefore $h^1({\cal K}^\vee\otimes{\cal
M}^{\otimes2})$ can be non-zero only if 
$$h^2({\cal K}^\vee\otimes{\cal
M}^{\otimes2})= h^0({\cal K}^{\otimes 2}\otimes{\cal
M}^{\otimes -2})\ne 0\ .$$
 This happens if and only if  ${\cal
M}^{\otimes 2} \simeq {\cal K}^{\otimes 2}$, i. e. ${\cal M}\otimes{\cal
K}^{-1} \in\Tors_2(\Pic)$. 

The surjectivity of $\Fg$ follows now from Proposition
\ref{filtrable}, and  the fact that trivial extensions cannot be
simple; the injectivity is stated in  2. and 3.

\qed

\begin{re}\label{us} The bundles ${\cal E}_{\cal L}$, ${\cal A}_{\cal R}$ admit a
unique rank 1 subsheaf with torsion free quotient (namely  ${\cal L}$ in the first
case and ${\cal K}\otimes{\cal R}$ in the second). 
\end{re}
\pf Indeed, if  ${\cal M}$ is   a
subsheaf with torsion free quotient of ${\cal E}_{\cal L}$ (${\cal A}_{\cal R}$) it
follows, as in the proof of Proposition
\ref{filtrable} that
$c_1^\Q({\cal M})\in\{0,c_1^\Q({\cal K})\}$.

If ${\cal M}$ was not ${\cal L}$ (respectively ${\cal K}\otimes{\cal R}$), one would
get a non-trivial morphism ${\cal M}\to {\cal K}\otimes{\cal L}^{-1}$ (respectively
${\cal M}\to{\cal R}$) which can be lifted to ${\cal E}_{\cal L}$ (respectively ${\cal
A}_{\cal R}$).  This would imply either the existence of curves in the forbidden
rational cohomology classes, or would give a splitting of the extension  which defines
${\cal E}_{\cal L}$ (${\cal A}_{\cal R}$).   
\qed 
\begin{co}\label{non-filtr} Let $\rho\in H^1(X,\Z_2)\setminus\{0\}$ and ${\cal L}_\rho$
the associated flat line bundle (see section \ref{roots}). Let ${\cal E}$ be a simple
bundle on $X$ with
$\det({\cal E})\simeq{\cal K}$, $c_2({\cal E})=0$ such that 
\begin{equation}\label{invo-inv}
{\cal E}\simeq {\cal E}\otimes{\cal L}_\rho\ .
\end{equation}
 Then ${\cal E}$ is
non-filtrable.
\end{co}
\pf  If ${\cal E}$ was filtrable  then, by Proposition \ref{bundles}, it would be
isomorphic with one of the bundles ${\cal E}_{\cal L}$, ${\cal A}_{\cal R}$. But, since
${\cal L}_\rho^{\otimes 2}\simeq {\cal O}_X$,  one has obviously
$${\cal E}_{\cal L}\otimes {\cal L}_\rho\simeq {\cal E}_{{\cal L}\otimes{\cal L}_\rho}\
,\ {\cal A}_{{\cal R}}\otimes{\cal L}_\rho\simeq {\cal A}_{{\cal R}\otimes{\cal L}_\rho}\
.
$$
Since ${\cal L}_\rho$ is non-trivial (see (\ref{iso})), it follows by Proposition
\ref{bundles}  that
${\cal E}_{\cal L}\not\simeq {\cal E}_{{\cal L}\otimes{\cal L}_\rho}$ and ${\cal A}_{{\cal
R}}\not\simeq{\cal A}_{{\cal R}\otimes{\cal L}_\rho}$. Therefore, the bundles ${\cal
E}_{\cal L}$, ${\cal A}_{\cal R}$ do not  verify  (\ref{invo-inv}), so they cannot by
isomorphic to ${\cal E}$.
\qed 
 
\subsection{Topological properties}
\begin{pr}
The isomorphism classes $[{\cal A}_{\cal R}]$ and $[{\cal
E}_{\cal R}]$ are not separable by disjoint neighborhoods in ${\cal M}^s$. 
More precisely, there exists an open neighborhood  $C'_{\cal R}$ of the
smooth point
$[{\cal A}_R]$ such that, with the notations of the previous proposition, one has
$$C'_{\cal R}\setminus\{[{\cal A}_{\cal R}]\}\subset C_{\cal R}\setminus\{[{\cal E}_{\cal
R}]\}\ .
$$
\end{pr} 
\pf  Let 
$$\delta_0=\left(\matrix{ \delta'_0&\sigma_0\cr 0&\delta''_0}\right)
$$
be an integrable semiconnection defining the holomorphic structure ${\cal E}_{{\cal
R}}$, as in the proof of Proposition \ref{bundles}, 5. Let  
$$\C\supset B\to{\cal A}^{0,1}_\kappa(\Theta\oplus [K\otimes\Theta])\ ,\
\delta_t=\left(\matrix{ \delta'_t&\sigma_t\cr \tau_t&\delta''_t}\right)$$
be a holomorphic map on the open disk such that $\delta(0)=\delta_0$, $\delta_t$ is
integrable,  and
$t\mapsto [\delta_t]$ parameterizes biholomorphically the curve $C_{\cal R}$. In
particular, the Dolbeault class $v:=[\dot\delta_0]\in H^1({\cal E}nd_0({\cal E}_{\cal
R}))=T_{[{\cal E}_{\cal R}]}({\cal M}^s)$ of the derivative
$\dot
\delta_0$ at 0 is a generator of the line $T_{[{\cal E}_{\cal R}]}(C_{\cal R})$.

The integrability condition implies
$$ \tau_t\circ \delta'_t+\delta''_t\circ
\tau_t=0\ .
$$
Differentiating at $0$ and taking into account that $\tau(0)=0$,
 we get
$$\dot\tau_0\circ \delta'_0+\delta''_0\circ\dot\tau_0=0\ .
$$
Since $\delta'_0\otimes\delta''_0=\kappa$, we see that
$\dot\tau_0$ is a
$\kappa$-closed
 $K$-valued $(0,1)$-form. Its Dolbeault cohomology class is obviously $h_*(v)$ (with the
notations in the proof of Proposition
\ref{bundles}, 5.), which does not vanish, because $C_{\cal R}$ is transversal to
$\Phi_{\cal R}$ in $[{\cal E}_{\cal R}]$.

For   $t\in B\setminus\{0\}$ set
$$g_t:=\left(\matrix{ t^{\frac{1}{2}}&0\cr 0&(t^{\frac{1}{2}})^{-1}}\right)\in
\qmod{SL(\Theta\oplus [K\otimes\Theta])}{\pm \rm Id}\ ,\ \delta^t:= g_t\cdot \delta_t
$$
One gets easily  
$$\lim_{t\ra 0} \delta^t= \left(\matrix{ \delta'_0&0\cr
\dot\tau_0&\delta''_0}\right)
$$
which defines precisely the holomorphic structure ${\cal A}_{\cal R}$, because the
Dolbeault class of $\dot \tau_0$ is non-zero. It suffices to notice that $\delta_t$
and $\delta^t$ define the same point in the moduli space, for every $t\ne 0$.
\qed\\
{\bf Remark:} An alternative proof can be obtained by studying the versal deformation
of the split bundle ${\cal R}\oplus[{\cal K}\otimes{\cal R}]$.

\begin{pr}\label{regular} The moduli space ${\cal M}^s$ is a smooth complex manifold of
dimension 1  at any point excepting the points
$[{\cal E}_{\cal R}]$, $[{\cal R}]\in\Tors_2(\Pic)$.
\end{pr}
\pf If ${\cal M}^s$ was not smooth of dimension 1 at $[{\cal E}]$, then
$H^0({\cal K}\otimes{\cal E}nd_0({\cal E}))\ne 0$. Let $\varphi\in H^0({\cal
K}\otimes{\cal E}nd_0({\cal E}))\setminus\{0\}$. Regarding $\varphi$ as a ${\cal
K}$-valued endomorphism, one obtains a section $\det(\varphi)\in H^0({\cal
K}^{\otimes 2})$, which must vanish, because  our surface has  $\kod(X)=-\infty$. 
Therefore, $\varphi$ has rank 1 which shows that ${\cal E}$  is filtrable. The
claim follows now from Proposition \ref{bundles}, 4., 5., 6.
\qed

\section{A moduli space of stable bundles}
\subsection{Classifying filtrable   stable bundles}

Let $g$ be a Gauduchon metric on $X$, and let ${\cal
M}^{\rm st}:={\cal M}^{\rm st}_g(0,{\cal K})$ be the moduli space   of
$g$-stable bundles with $c_2=0$ and $\det={\cal K}$. 

We denote by $\Pic^c_{<t}$ ($\Pic^T_{<t}$) the subset of
$\Pic^c$ (respectively $\Pic^T$)  defined by the inequality $\deg_g({\cal
L})<t$. Each $\Pic^c_{<t}$ is an open {\it pointed disk}, i.e. an open  disk
minus a   point; $\Pic^T_{<t}$ is a finite union of pointed disks. We define
similarly the spaces 
$$\Pic^c_{\leq t}\ ,\ \Pic^T_{\leq t}\ ,\ \Pic^c_{=t}\ ,\ \Pic^T_{=t}\ ,\  
\left[{\Pic^c}\right]_{\leq t}^{\geq s}\ ,\  \left[{\Pic^T}\right]_{\leq t}^{\geq s}\
.$$
\vspace{2mm}\\
{\bf Assumption:} {\it We   assume that the metric $g$ was chosen such that
$\deg_g({\cal K})\ne 0$. }\vspace{1mm}

The set of Gauduchon metrics satisfying this assumption is open and
dense.   Indeed, since  $c_1({\cal K})^2\ne 0$ , it follows
that,  perturbing 
$\omega_g$ with a  generic,  closed, cohomologically non-trivial, small 
$(1,1)$-form, will yield such a metric. The same argument applies to any
given holomorphic line bundle ${\cal L}$ with non-trivial  $c_1({\cal L})^2$.
Note however, that, in general, it  is very difficult to see whether $\deg_g({\cal
L})$  is always positive, always negative, or can have both signs as $g$ varies in the
space of Gauduchon metrics. For instance, if the $\partial\bar\partial$-Chern class of
${\cal L}$ is a rational combination with positive coefficients of
$\partial\bar\partial$-Chern classes of line bundles associated with   curves, the
degree
$\deg_g({\cal L})$  will be always positive. 

On the other hand -- as noticed by one of the two referees  -- using Buchdahl's   
ampleness criterion for non-K\"ahlerian surfaces \cite{Bu2}, one  can prove that   both signes
are possible when
$c_1({\cal L})^2\ne 0$ and $X$ has no curves. Therefore, under this  
hypothesis on $X$, one can choose $g$ such that $\deg_g({\cal
K})<0$ and continue the proof   treating only this case. Since  the case
$\deg_g({\cal K})>0$ is not much more difficult than the other one,  we will
not follow this way.\\

Set $\kg:=\deg_g({\cal K})/2$. Remark \ref{us} shows that the stability of ${\cal
E}_{\cal L}$ reduces to the condition $\deg_g({\cal L})<\kg$, and the stability of 
${\cal A}_{\cal R}$ reduces to the condition $\deg_g({\cal K})<\kg$. 
Therefore

\begin{thry}\label{st} Under the assumptions and with the notations of Proposition
\ref{bundles} the following holds:
 
\begin{enumerate}
\item ${\cal E}_{\cal L}$ is $g$-stable if and only if $\deg_g({\cal L})<\kg$.
\item When $\deg_g{\cal K}<0$, the bundles ${\cal A}_{\cal R}$, $[{\cal
R}]\in\Tors_2(\Pic)$ are all stable. When $\deg_g{\cal K}>0$, they are not stable.
\item If $\deg_g({\cal K})<0$, then the restriction
$$\resto{\Fg}{\Pic^T_{<\kg}\coprod\Tors_2(\Pic)}: \Pic^T_{<\kg}\coprod
 \Tors_2(\Pic)\map {\cal M}^s
$$
maps bijectively $\Pic^T_{<\kg}\coprod
 \Tors_2(\Pic)$ on the filtrable part of ${\cal M}^{\rm
st}$. The image of the subspace $\Pic^T_{<\kg}$ is open in ${\cal M}^{\rm
st}$.
\item If $\deg_g({\cal K})>0$, then the restriction
$$\resto{\Fg}{\Pic^T_{<\kg}}: \Pic^T_{<\kg}\map {\cal M}^s
$$
maps bijectively $\Pic^T_{<\kg}$ on the filtrable part of ${\cal M}^{\rm
st}$. The image of the subspace $\Pic^T_{<\kg}\setminus\Tors_2(\Pic)$ is open in ${\cal
M}^{\rm st}$.

\end{enumerate}
\end{thry}

Throughout the rest of the paper $X$ will always denote a class $VII$ surface with
$b_2=1$ satisfying the hypothesis of  Proposition
\ref{bundles}.

\def\ASD{{\rm ASD}}
 
\subsection{A collar around the reductions}

Let $E$ be a Hermitian 2-bundle with $c_2(E)=0$, $\det(E)=K$, and let $a$ be the Chern
connection of the pair $({\cal K},\det(h))$.

For every $c\in\Tors(H^2(X,\Z))$ one has a unique (up to the gauge group
$SU(E)$) orthogonal decomposition
$$E=L\oplus [K\otimes L^\vee]\ ,
$$
where $L$ is a Hermitian line bundle of Chern class $c$, and these are the only
splittings of $E$ (by the same
computation as in the proof of Proposition
\ref{filtrable}).   Every  reducible connection 
  $A\in {\cal M}^{\ASD}_a(E)$ is equivalent to a direct sum $b\oplus (a\otimes
b^\vee)$, where $b$ is a Hermitian connection on $L$ such that 
\begin{equation}\label{asd}
(2F_b-F_a)^+=0\ .
\end{equation}
The condition $(2F_b-F_a)^+=0$  means that $2F_b-F_a$ is an ASD form  so, being closed,
it coincides with the unique harmonic representative $r$ of the de Rham cohomology class
$-2\pi ic_1^{\rm DR}(X)$. The equation (\ref{asd}) becomes
\begin{equation}\label{ab}
F_b=\frac{1}{2}(F_a+r)\ .
\end{equation}
Let ${\cal N}_c$ be the moduli space of solutions $b$ of (\ref{ab}) modulo the gauge
group ${\cal C}^\infty(X,S^1)$ of $L$. It is well-known that the moduli space of
Hermitian connections of fixed curvature on a Hermitian line bundle over an arbitrary 
compact manifold $V$ is a $\T_V$-torsor\footnote{In general, a $G$-torsor
is a set $\Gamma$ endowed with a free transitive $G$-action. Fixing a point
$\gamma\in\Gamma$ gives an identification $\Gamma\simeq G$, and any two such
identifications differ by a $G$-translation.}, where
$\T_V$ is the torus
$$\T_V:=\qmod{iH^1(V,\R)}{2\pi iH^1(V,\Z)}\ .
$$
Therefore it is (non-canonically) isomorphic to this torus. In our case,
$\T_X$ is obviously a circle. Therefore
\begin{pr}  The subspace ${\cal M}^{\rm red}\subset {\cal M}^{\ASD}_a(E)$  of
reducible solutions decomposes as a disjoint union 
$${\cal M}^{\rm red}=\coprod_{c\in \Tors(H^2(X,\Z))}{\cal M}^{\rm red}_c
$$
where ${\cal M}^{\rm red}_c$ is the moduli space of solutions which admit a
Hermitian line bundle   $L$ with $c_1(L)=c$ as  parallel summand. 
Each ${\cal M}^{\rm red}_c$ is naturally isomorphic to the circle ${\cal N}_c$. 
\end{pr}
Our next purpose is to understand the topology of ${\cal M}^\ASD_a(E)$ around the
reducible locus ${\cal M}^{\rm red}$. It is very natural (and almost obvious) that the
circle ${\cal N}_c$ is contained in the closure of the pointed disk
$\Fg(\Pic^c_{<\kg})$, so that the union  $\Fg(\Pic^c_{<\kg})\cup {\cal N}_c$
becomes a closed pointed  disk. 
Indeed, the limit of the stable bundle ${\cal E}_{\cal
L}$ as ${\cal L}\ra {\cal L}_0$ (with $\deg_g({\cal L})_0=\kg$) should be the
polystable bundle  ${\cal L}_0\oplus [{\cal K} \otimes{\cal L}^\vee_0]$, and the circle
of polystable bundles of this form can be naturally identified with ${\cal N}_c$.   

However,   the local structure of a moduli space of polystable bundles around
the non-stable points can be in general very complicated, and it is not a
subject available in the literature for our non-K\"ahlerian framework. Hence we
will indicate a simple ad-hoc argument. 
\begin{lm}\label{smooth} Let $A$ be a reducible solution in $	{\cal A}^{\ASD}_a(E)$.
Then the harmonic space $\H^2_A$ of the deformation elliptic  complex associated with
$A$  vanishes.
\end{lm}
\pf  Let $L\oplus [K\otimes L^\vee]$
be the decomposition of $E$ in
$A$-parallel factors, and let $A=b\oplus c$ be the corresponding
splitting of $A$. One has an 
$A$-parallel splitting 
$${\rm su}(E)=(X\times[i\R])\oplus [K^\vee\otimes L^2]
$$
and the induced connections on the summands are the trivial connection and 
 $f:=b\otimes c^\vee= a^\vee\otimes b^{\otimes 2}$ respectively. The connections
$b$ and $f$ are integrable, so they define holomorphic structures ${\cal L}$, ${\cal
F}$ on
$L$ and
 $F:= K^\vee\otimes L^2$. Note   $f$ is an ASD connection (or, equivalently, a
Hermitian-Einstein connection of vanishing Einstein constant), and that
$$\deg_g({\cal L})=\kg\ ,\ \deg_g({\cal F})=0\ .
$$

As in the K\"ahlerian case \cite{K}, the idea is to compare the
$d^+$-complex of
$f$ with the Dolbeault complex of $\bar\partial_f$.   We have the diagram
$$
\begin{array}{ccccccccc}
0&\ra&0&\map&A^{10}(F)&\textmap{(\partial_f,p^\omega\bar\partial_f) }&A^{20}(F)\oplus
A^0(F)\omega_g&\ra& 0
\\ \\
&&\downarrow&&\downarrow i&&\downarrow j
\\ \\
0&\ra&A^0(F)&\textmap{d_f}&A^1(F)&\textmap{d_f^+}&A^2_+(F)&\ra& 0
\\ \\
&&\parallel&&\downarrow p^{01}&&\downarrow p^{02}
\\ \\
0&\ra&A^0(F)&\textmap{\bar\partial_f}&A^{01}(F)&\textmap{\bar\partial_f}&A^{02}
(F)&\ra& 0
\end{array}
$$
where $i$, $j$ are the obvious inclusions, and  $p^{01}$, $p^{02}$, $p^\omega$ the
obvious projections.

The cohomology of the third line is just the cohomology of the holo\-morphic bundle
${\cal F}$ associated with the semiconnection $\bar\partial_f$. But
$$H^2({\cal F})=H^2({\cal K}^\vee\otimes{\cal L}^{\otimes 2})=H^0({\cal K}^{\otimes
2}\otimes{\cal L}^{\otimes -2})^\vee=0
$$
because our surface does not have numerically bicanonical curves. It suffices to
show that the second cohomology of the first line vanishes.

Let $(a^{20},\varphi\omega_g)\in\coker(\partial_f,p^\omega\bar\partial_f)$. This
implies
$$\partial_f^*(a^{20})+\bar\partial_f^*(\varphi\omega_g)=0\ ,
$$
hence, since $a^{20}$ and $\omega_g$ are self-dual with respect to the ($\C$-linear)
Hodge operator, 
$$\bar\partial_f(a^{20})+ \partial_f (\varphi \omega_g)=0\ .
$$
We get \footnote{I am grateful to one of the two referees for pointing me out an error
(caused by a missing term) in this part of the proof, and for indicating a correct
argument.}
$\bar\partial_f\partial_f(
\varphi\omega_g)=0$. We claim that the kernel of the second order elliptic   operator 
$q_f:A^0(F)\to  A^{0}(F)$ defined by
 $q_f(\psi):=i*\bar\partial_f\partial_f(
\psi\omega_g)$ 
 vanishes.  The adjoint of this operator is  
$p_f=i\Lambda_g \bar\partial_f\partial_f$ (see \cite{Bu1}, \cite{LT1} p. 225).  Taking
into account $i\Lambda_g(\bar\partial_f\partial_f+\partial_f\bar\partial_f)=0$ (because
$f$ is ASD), one gets the identity
$$i\Lambda_g\bar\partial\partial
|\psi|^2=(p_f(\psi),\psi)-|d_f(\psi)|^2+
(\psi,p_f(\psi)) 
$$
(compare with the computation in \cite{LT1}, Lemma 1.2.5, p. 31). By the maximum principle,
it follows that $\ker(p_f)=\ker (d_f)$, which vanishes,  because ${\cal F}$ is
non-trivial. Since   both operators have index 0, we also get $\ker(q_f)=0$. Therefore
$\varphi=0$, which implies that $\bar\partial_f(a^{20})=0$, so $a^{20}$ defines a
holomorphic section in
${\cal K}\otimes {\cal F}={\cal L}^{ \otimes 2}$. It follows that
$a^{20}=0$, because the holomorphic line bundle ${\cal L}^{ \otimes 2}$
is non-trivial (has non-vanishing degree) and our surface has no
numerically trivial curve.
\qed

For every holomorphic line bundle ${\cal L}$ we denote by $b_{\cal L}$ its unique
Hermitian-Einstein connection.

\begin{pr}\label{collar}

 The extension
$\bar\Fg:\Pic^T_{\leq \kg}\to {\cal M}^{\ASD}_a(E)$
of $\resto{\Fg}{\Pic^T_{<\kg}}$ defined by
$$\bar\Fg([{\cal L}]):=[b_{\cal L}\oplus (a\otimes b_{\cal L}^\vee)]\ ,\ \forall
[{\cal L}]\in\Pic_{=\kg}^T
$$
maps homeomorphically $\Pic^T_{\leq \kg}\setminus\Tors_2(\Pic)$ on an open subspace
of ${\cal M}^{\ASD}_a(E)$. In particular  ${\cal M}^{\ASD}_a(E)$ has the structure of
a real 2-manifold with boundary  near ${\cal M}^{\rm red}$.
\end{pr}
\pf  By the standard Kobayashi-Hitchin correspondence between stable bundles and
Hermitian-Einstein connections \cite{LT1} we get immediately  that the restriction
$\resto{\Fg}{\Pic^T_{<\kg}\setminus\Tors_2(\Pic)}$ is an open embedding.

It remains to prove that $\bar\Fg$ is a local homeomorphism near the union of circles
$\Pic_{=\kg}^T$.

Let $\varepsilon >0$.  
For an integrable connection $A$ denote by ${\cal E}_A$ the holomorphic bundle
defined by the semiconnection $\bar\partial_A$. We define
$${\cal M}^{\ASD}_a(E)_{\varepsilon}:=
 \{[A]\in {\cal M}^{\ASD}_a(E)|\ \exists
[{\cal L}]\in[\Pic^T]_{\leq \kg}^{\geq\kg-\varepsilon}\hbox{ with } H^0({\cal
L}^\vee\otimes{\cal E}_A)\ne 0\}\ .
$$

By elliptic semicontinuity it follows easily that ${\cal M}^{\ASD}_a(E)_{
\varepsilon}$ is a closed (hence  compact by Corollary \ref{compactness}) subspace of 
${\cal M}^{\ASD}_a(E)$. 

The map
$$\lambda_\varepsilon:{\cal M}^{\ASD}_a(E)_{\varepsilon}\map
 \left[{\Pic^T}\right]_{\leq
\kg}^{\geq\kg-\varepsilon} 
$$
defined by
$$[A]\mapsto \hbox{the unique }[{\cal L}]\in\Pic^T\hbox{ such that }H^0({\cal
L}^\vee\otimes{\cal E}_A)\ne 0
$$
is continuous, by elliptic semicontinuity again, and bijective. The uniqueness of ${\cal
L}$ follows easily from our hypothesis (non-existence of numerically trivial curves).
Therefore, by   an elementary topological lemma, this map is a homeomorphism. By
definition, this homeomorphism is the inverse of the restriction
 $\resto{\bar\Fg}{\left[\Pic^T\right]_{\leq
\kg}^{\geq \kg-\varepsilon}}$. It suffices to prove that
\\ \\
{\bf Claim:} ${\cal M}^{\ASD}_a(E)_{\varepsilon}$ is a neighborhood of the reducible
part ${\cal M}^{\rm red}$.
\\

Let $A:=b\oplus c$ be a reducible connection in ${\cal A}^{\rm ASD}_a(E)$ as in the
proof of Lemma \ref{smooth},
$E=L\oplus [K\otimes L^\vee]$ be the corresponding $A$-parallel decomposition,
$F:=K^\vee\otimes L^{\otimes 2}$,  $f:= a^\vee\otimes b^{\otimes 2}$, and ${\cal F}$
the  corresponding holomorphic line bundle.

Consider the   space  ${\cal S}$ of solutions of the following elliptic differential
system
\begin{equation}\label{s}
\left\{
\begin{array}{ccc}
\bar\partial\beta^{01}&=&0\\
\Lambda_g(d\beta -\alpha \wedge \bar\alpha )&=&0\\
\bar\partial_{f} \alpha +2\beta^{01}\wedge\alpha&=&0\\
\Lambda_g(\partial_{f}\alpha +2\beta^{10}\wedge \alpha)&=&0\\
d^* \beta&=&0\\
\end{array}
\right.\ ,\ \beta\in iA^1(X)\ ,\ \alpha\in A^{01}(F)\ .
\end{equation}

This space  comes with a natural $S^1$-action:
$$(\zeta,(\beta,\alpha))\mapsto (\beta,\zeta^2 \alpha) \ .
$$

The linearization of   (\ref{s}) at $0$ is  
\begin{equation} 
\left\{
\begin{array}{ccccccc}
d^+\beta &=&0&,&d^*\beta&=&0\\
 \bar\partial_{f} \alpha  &=&0&,&\Lambda_g \partial_{f}\alpha &=&0\\ 
\end{array}
\right.\    ,
\end{equation}
so the tangent space $T_0({\cal S})$   can be identified with  $i H^1(X,\R)\oplus
H^{1}({\cal F})$. Indeed, by the maximum principle it is easy to prove that any
$(0,1)$-Dolbeault cohomology class of
${\cal F}$ has a unique representative $\alpha$ with
$\Lambda_g\partial_{f}\alpha=0$ (see the proof of  Lemma \ref{smooth}).
The map defined by the  left hand side of the system (\ref{s}) is a
submersion at 0, hence ${\cal S}$ is a smooth manifold of dimension 3
around 0.

We define the map $\Ag:{\cal S}\ra {\cal A}_a^{\ASD}(E)$ by
$$\Ag(\beta,\alpha):=\left(\matrix{b+\beta&\alpha\cr-\bar\alpha&c-\beta}\right)\ .
$$
We consider ${\cal A}_a^{\ASD}(E)$ as an $S^1$-space via $\zeta\mapsto
\left(\matrix{\zeta&0\cr0& \zeta^{-1}}\right)\in\SU(E)$. Note that, by Lemma
\ref{smooth},
${\cal A}_a^{\ASD}(E)$ is a smooth Banach manifold at $A$ (after suitable Sobolev
completion). It is easy to check that 
\begin{enumerate}
\item $\Ag$ is $S^1$ equivariant and $\Ag(0)=A$,
\item $\Ag$ is an immersion at $0$,
\item $\Ag_*(T_0({\cal S}))$ is a complement of $T_A[A]$ in $T_A({\cal
A}_a^{\ASD}(E))$; in particular 
$\im(\Ag)$ is transversal in
$A$ at the
$\SU(E)$-orbit
$[A]=\SU(E)\cdot A$.
\end{enumerate}

It follows that $\im(\Ag)$ is mapped on a neighborhood of   $[A]$ in ${\cal
M}^{\ASD}_a(E)$ and that ${\cal S}/S^1$ is a local model for the moduli space
 ${\cal M}^{\ASD}_a(E)$ at $[A]$.
Therefore, in order to prove our claim, it suffices to show that, when
$(\beta,\alpha)\in{\cal S}$ are sufficiently small, then $[\Ag(\beta,\alpha)]$ belongs
to ${\cal M}^{\ASD}_a(E)_\varepsilon$.

But, by construction, since $\bar\alpha$ is of type $(1,0)$, the summand $L$ of $E$ is
always a holomorphic sub line bundle of the holomorphic bundle ${\cal
E}_{\Ag(\beta,\alpha)}$. The induced holomorphic structure on this summand is defined
by the semiconnection
$\bar\partial_b+\beta^{01}$. The second equation in the system (\ref{s}) shows that
the degree of this holomorphic structure is just $\|\alpha\|^2_{L^2}$. Therefore, this
degree becomes smaller than
$\varepsilon$ for $\alpha$ sufficiently small.
\qed

\subsection{The missing centers}

\def\Pgic{{\Pg}{\rm ic}}
Denote by $\Pic^c_-$, respectively $\Pic^c_+$, the   space (isomorphic to $\C$) obtained
by adding formally a point
$0_c$ (respectively $\infty_c$) to $\Pic^c$ such that the limits 
$$\lim_{\begin{array}{cc}\scriptstyle c_1({\cal L})=c\vspace{-1mm}\\
\scriptstyle\deg_g({\cal L})\to-\infty\end{array}}[{\cal
L}]=0_c\ ,\ \lim_{\begin{array}{cc}\scriptstyle c_1({\cal L})=c\vspace{-1mm}\\
\scriptstyle\deg_g({\cal L})\to+\infty\end{array}}[{\cal
L}]=\infty_c  
$$
exist.  Set
$$ \Pic^T_\pm:=\union_{c\in
\Tors(H^2(X,\Z))}\Pic^c_\pm\ .
$$
As in the previous section, we
introduce the subspaces
$$[\Pic^c_-]_{< t} ,\ [\Pic^T_-]_{< t}  , [\Pic^c_-]_{\leq t} ,\ [\Pic^T_-]_{\leq t} 
,\ [\Pic^c_+]_{>t}  ,\ [\Pic^T_+]_{> t}  , [\Pic^c_+]_{\geq t}  ,\ [\Pic^T_+]_{\geq
t}\ .
   $$
\begin{pr} 
\begin{enumerate}\label{extension-map}
\item  ${\cal M}^{\st}$ is a complex space of pure dimension 1, which is
smooth in the case $\deg_g({\cal K})<0$ and whose singular locus is the finite
set $\{[{\cal E}_{\cal R}]\}_{[R]\in\Tors_2(\Pic)}$ in the case $\deg_g({\cal K})>0$.

\item The map
$$ {\bar\Fg} :\Pic^T_{\leq \kg}\map {\cal
M}^\ASD_a(E)=\overline{{\cal M}^{\st}}
$$
extends to a map ${\tilde \Fg} :[\Pic^T_-]_{\leq \kg}\map {\cal
M}^\ASD_a(E)=\overline{{\cal M}^{\st}}$ which is holomorphic on $[\Pic^T_-]_{<\kg}$.

\item The points $\tilde\Fg(0_c)$ are fixed under the  action of  
$H^1(X,\Z)/2H^1(X,\Z)\simeq\Z_2$ on the moduli space (see sections \ref{donaldson},
\ref{roots}).
\end{enumerate}
\end{pr}
\pf 1. The first statement follows directly from Proposition \ref{regular} and
Proposition
\ref{bundles}.

2. By Corollary \ref{compactness}, the space ${\cal M}^\ASD$ obtained by adding to
${\cal M}^\st$ the circles of polystable bundles is compact.  The idea of the proof is
based on the following simple fact: the only way to compactify a complex line
$\P^1(\C)\setminus\{0\}$ {\it as a 1-dimensional complex space} is   by adding a point
(which a priori could be of course singular).

Consider the  compact 1-dimensional complex space 
$${\cal C}:=[\Pic^T_+]_{> \kg-\varepsilon}\ {\union}_{\Fg} {\cal M}^\st
$$
obtained by filling in disks in the direction $\deg_g({\cal L})\to\infty$. By
Proposition
\ref{collar}, and the first statement in this proposition, ${\cal C}$ is a complex
space of dimension 1 with at most $|\Tors_2(\Pic)|$ simple (normal crossing)
singularities.

The irreducible component of ${\cal C}$ which contains the line 
$$[\Pic^c_+]_{>
\kg-\varepsilon}\ {\union}_{\Fg}\bar\Fg(\Pic^c_{<\kg})$$
 is an irreducible algebraic
curve, so it can be embedded in a projective space $\P^n(\C)$ such that the
intersection with the hyperplane at infinity is 
$\{\infty_c\}$. Then $\resto{\Fg}{{\Pic^c_{<\kg}}}$  factorizes through a bounded 
$\C^n$-valued map, which can be obviously extended holomorphically in $0_c$.

3. Let $\otimes\rho:{\cal M}^\st\to {\cal M}^\st$ be the involution induced by the
generator $\rho$ of $H^1(X,\Z)/2H^1(X,\Z)$, and $[{\cal L}_\rho]\in\Pic^0$ be the
line bundle associated with the representation $\pi_1(X)\ra\{\pm1\}\subset
S^1\subset\C^*$ defined by $\rho$ (see section \ref{roots}).   One has obviously the
identity
$$\otimes\rho(\Fg([{\cal L}])=[{\cal E}_{\cal L}\otimes{\cal L}_\rho]=\Fg([{\cal
L}\otimes{\cal L}_\rho])\ .
$$
It suffices now to note that 
$$\lim_{\begin{array}{cc}\scriptstyle c_1({\cal L})=c\vspace{-1mm}\\
\scriptstyle\deg_g({\cal L})\to-\infty\end{array}}[{\cal
L}\otimes{\cal L}_\rho]=0_c\ ,
$$
because $\deg_g({\cal L}_\rho)=0$.
\qed

The fact that the points $\tilde\Fg(0_c)$ are fixed under the involution $\otimes\rho$
is very important.  Using Corollary \ref{non-filtr} and Proposition \ref{regular}, one
obtains
\begin{co}\label{nonf} The point $\tilde\Fg(0_c)$ corresponds to a non-filtrable bundle,
  in particular  it is a smooth point in the
moduli space ${\cal M}^\st$.
\end{co}
Let ${\cal E}_c$ be a bundle in the isomorphy class $\tilde\Fg(0_c)$. The isomorphism
${\cal E}_c\simeq {\cal E}_c\otimes{\cal L}_\rho$ shows that ${\cal E}_c$ can be obtained
as push-forward of a line bundle on a bicovering. More precisely, let $\pi_\rho:\tilde
X_\rho\to X$ be the bicovering of $X$ associated with $\rho$ (regarded as a representation
$\pi_1(X,x_0)\to\Z_2$).  
\begin{re} For every $c\in \Tors(H^2(X,\Z))$ there exists a holomorphic line bundle
${\cal M}_c$ on  $\tilde X_\rho$ such that ${\cal E}_c\simeq  [\pi_\rho]_*({\cal
M}_c)$.
\end{re}
{\bf Remark:} Our results so far show that  the moduli space $\overline{{\cal
M}^\st}$ contains as an open set the disjoint union  
$\left[\coprod_{i=1}^{p-q} D_i\right]\coprod \left[\coprod_{j=1}^{q} P_j\right]
$
described in section \ref{strategy}. The point  $\tilde\Fg(0_c)$ corresponding to
a class $c\in\Tors_2(H^2(X,\Z))$ is the center  $O$ of the corresponding
``horizontal" disk
$\tilde\Fg([\Pic^c_-]_{\leq
\kg})$ (denoted by $E$ in the picture).  Corollary
\ref{nonf} plays an important role: it shows in particular that  the center of a
``vertical" disk  $C'_i$ (or $C_i$) cannot coincide with a center $\tilde\Fg(0_c)$, hence
such a vertical disk cannot be contained in the ``horizontal" locus 
$\im(\tilde\Fg)$.
\subsection{A smooth compact complex curve in the moduli space}
We can prove now:
\begin{thry}\label{curve} Let $X$ be a class $VII$ surface with $b_2=1$ which does not
admit any   divisor $C>0$ with $c_1^\Q({\cal O}(C))\in\{0,\pm c_1^\Q({\cal K}),
2c_1^\Q({\cal K)}\}$. Let $g$ be a Gauduchon metric on $X$ such that $\deg_g({\cal K})\ne
0$.
\begin{enumerate}
\item Suppose $\deg_g({\cal K})<0$.
Then 
$\overline{{\cal M}^\st}\setminus \tilde\Fg([\Pic^T_-]_{\leq\kg})$
is a \ub{non-empty}, possibly non-connected, smooth, closed complex curve whose
only filtrable points  are the extensions
$[{\cal A}_R]$, $[{\cal R}]\in\Tors_2(\Pic)$. 

\item Suppose $\deg_g({\cal K})>0$. The closure of $\overline{{\cal
M}^\st}\setminus \tilde\Fg([\Pic^T_-]_{\leq\kg})$ is a \ub{non-empty}, possibly
non-connected, smooth, closed complex curve whose only filtrable points 	are
the extensions
$[{\cal E}_R]$, $[{\cal R}]\in\Tors_2(\Pic)$.
 \end{enumerate}
\end{thry}
\pf  1. In this case, by Theorem \ref{st} and Proposition \ref{collar} the subspace
$\tilde\Fg([\Pic^T_-]_{\leq\kg})$ is open and closed in
$\overline{{\cal M}^\st}$. It suffices to note that
$[{\cal A}_{\cal R}]$ does not belong to this subspace. This follows from the
injectivity of the map $\Fg$ and Corollary \ref{nonf} which assures that the centers
$\tilde\Fg(0_c)$ are non-filtrable. 
\\ \\
2. For  any $[{\cal R}]\in\Tors_2(\Pic)$   the branch $C_{\cal R}$ passing through
$[{\cal E}_{\cal R}]$ (see  Proposition \ref{bundles}) consists of stable bundles,
because stability is an open property \cite{LT1}. It suffices to show that the
irreducible component containing this branch is not contained in
$\tilde\Fg([\Pic^T_-]_{\leq\kg})$. The normal crossing $C_{\cal R}\cup\Phi_{\cal R}$ is
mapped injectively in the moduli space. If $C_{\cal R}$ was contained in
$\tilde\Fg([\Pic^T_-]_{\leq\kg})$, the center
$[{\cal E}_{\cal R}]$ of this branch would coincide with an element of
$\tilde\Fg([\Pic^T_-]_{\leq\kg})\setminus\Phi_{\cal R}$, which is impossible 
(use again the  injectivity of the map $\Fg$ and Corollary \ref{nonf}). 
\qed
\\ \\
Theorem \ref{curve} yields a closed   curve $Y$ and a holomorphic morphism $Y\to {\cal M}^{\rm
st}$ taking both filtrable and non-filtrable values. In the next section we will see that  such a
morphism cannot exist.

\section{ Families of  bundles   parameterized by a
curve}\label{families}

We begin with the following result concerning the correspondence between holomorphic
morphisms from a curve into a moduli space of simple bundles    and
holomorphic families.

\begin{lm}\label{family} Let $X$ be a complex manifold, $E$ a differentiable rank 2
bundle on
$X$, and ${\cal L}$ a  fixed holomorphic structure on $\det(E)$. Let ${\cal
M}^\s(E,{\cal L})$ be the moduli space  of simple holomorphic
structures on
$E$ which induce ${\cal L}$ on $\det(E)$, modulo the complex gauge group ${\rm SL}(E)$.

Let $Y$ be a compact complex curve and $\fg:Y\to {\cal M}^\s(E,{\cal L})$ be a holomorphic
morphism. There exists a line bundle ${\cal N}$ on $Y$ and a holomorphic 2-bundle
${\cal E}$ on
$Y\times X$ such that
\begin{enumerate}
\item The family ${\cal E}$ induces $\fg$, i.e.  
$$\left[\resto{{\cal E}}{\{y\}\times X}\right]=\fg(y)\ ,\ \forall y\in Y\ .$$
\item $\det({\cal E})\simeq p_Y^*({\cal N})\otimes p_X^*({\cal L})$, where $p_X$, $p_Y$
are the respective projections.
\end{enumerate}
\end{lm}
\pf
 The   deformation theory for holomorphic bundles   extends easily to
holomorphic structures with fixed determinant.  In particular, the germ of the moduli
space
${\cal M}^\s(E,{\cal L})$ at a point
$[{\cal E}]$ is the basis of a universal deformation of ${\cal E}$ in the space of
holomorphic structures compatible with ${\cal L}$  (see \cite{Miy} for the case of
plain bundles).    This provides   
   an open cover
${\cal U}=(U_i)_{i\in I}$ of
$Y$ and bundles ${\cal E}_i$ on $U_i\times X$ with isomorphisms $j_i:\det({\cal
E}_i)\simeq [p_X^i]^*({\cal L})$, where $p_X^i:U_i\times X\to X$ are the projections on
the $X$-factor. 

If the cover ${\cal U}$ is sufficiently fine, on the intersections $U_i\cap U_j$ we get
$SL$-isomorphisms
$f_{ji}:{\cal E}_i\to {\cal E}_j$. 
These isomorphisms are obtained in the usual way, by noting that
  the sheaves $(p_{U_i\cap U_j})_*({\cal E}_i^\vee\otimes{\cal E}_j)$ are
line bundles, which will be trivial if the cover is sufficiently fine. If the
intersections $U_i\cap U_j$ are simply connected one can find sections in these
trivial line bundles which correspond to $SL$-isomorphisms.

The compositions  $f_{kji}=f_{kj}\circ f_{ji} \circ
f_{ki}^{-1}$ form a 2-cocycle with coefficients in $\Z_2$. Its cohomology class  
$w(\fg,{\cal E}_i,f_i)\in\check{H}^2({\cal U},\Z_2)$ is the obstruction to gluing
 the pairs
$({\cal E}_i,j_i)$ to   a global family ${\cal E}$ on $Y\times X$ with a global
isomorphism
$j:\det({\cal E})\simeq p_X^*({\cal L})$. The corresponding class   
$w(\fg)\in H^2(Y,\Z_2)$ depends only on $\fg$ and measures the obstruction to the
existence of a family ${\cal E}$ on $Y\times X$ with an isomorphism
$j:\det({\cal E})\simeq p_X^*({\cal L})$ inducing $\fg$. 

Since $H^2(Y,{\cal O}_Y^*)=0$, after passing to a finer cover if necessary,  we can
find a \v{C}ech 1-cochain $\eta=(\eta_{ij})_{i,j} \in
\check{C}^1({\cal U},{\cal O}_Y^*)$ such that $\delta (\eta)=w(\fg,{\cal
E}_i,f_i)$. The system $(\eta_{ji}f_{ji} )_{j,i}$ satisfies the cocycle
condition, so we can glue  the bundles ${\cal E}_i$ via this system of
isomorphisms and get a global bundle ${\cal E}$ on $Y\times X$ inducing the map
$\fg$. The determinant of this family will be $p_Y^*({\cal N})\otimes p_X^*({\cal
L})$, where ${\cal N}$ is the line bundle associated with the cocycle
$(\eta_{ij}^2)_{i,j}$. 
\qed
 
\begin{thry}\label{constant} Let $X$ be  a surface of algebraic dimension $a(X)=0$, and ${\cal E}$ a holomorphic 2-bundle    on
$Y\times X$. There
exists a non-empty Zariski open set
$U\subset X$, a coherent sheaf  ${\cal T}$ of rank 1 or 2 on $X$ which is locally free on $U$ 
and, for every
$y\in Y$,   a morphism
$e_y:{\cal T}\to {\cal E}_y$ which is a bundle embedding (i.e. fibrewise injective) on $U$.  
\end{thry}
\pf\footnote{The short proof we reproduce here was kindly suggested by one of the referees.}  For  
any holomorphic 2-bundle
${\cal F}$ on the curve
$Y$ we put
$$d({\cal F}):=\min\left\{d\in\Z |\ \exists [{\cal M}]\in\Pic^d(Y) \hbox{ s. t. } H^0({\cal M}\otimes
{\cal F})\ne 0\right\}\ , $$
$$z({\cal F}):=\left\{[{\cal M}]\in\Pic^{d({\cal F})}(Y)|\ H^0({\cal
M} \otimes {\cal F})\ne 0\right\}\subset \Pic^{d({\cal F})}(Y)\ .
$$

Note first that   every $[{\cal M}]\in z({\cal F})$   has the following two important properties
\begin{equation}\label{P}  
\left\{\begin{array}{c} \ker\left[ {\rm ev}_y:H^0({\cal M} \otimes
{\cal F}) 
\to {\cal M}(y)\otimes{\cal F}(y))\right]= 0\ \forall y\in Y,   \vspace{2mm} \\ 
h^0({\cal M}\otimes {\cal
F})\in\{1,2\}\ .
\end{array}\right.
\end{equation}
 Indeed, if   $s\in H^0({\cal M} \otimes {\cal F})\setminus\{0\}$ vanished at $y$,  it would
induce a nontrivial section in $H^0({\cal M}(-y) \otimes {\cal F})$, contradicting the minimality of
$d({\cal F})$. The second formula follows from the first and the definition of $z({\cal F})$.
\\

 For   $x\in X$ denote ${\cal E}^x:=\resto{{\cal E}}{Y\times \{x\}}$
(regarded as a bundle on $Y$). The sets 
$$X_d:=\{x\in X|\ d({\cal E}^x)\leq d\}
$$
are Zariski closed in $X$ and the map $x\mapsto d({\cal E}^x)$ is bounded from above. Put
$$\delta:=\max_{x\in X} d({\cal E}^x) =\min\{d\in\Z| \forall x\in X\ \exists [{\cal
M}]\in\Pic^d(Y) \hbox{ s. t. } H^0({\cal M} \otimes {\cal E}^x)\ne 0\}\ ,
$$
$$W:=\{x\in X|\ d({\cal E}^x)=\delta\}\ ,
$$
$$Z:=\{(x,[{\cal M}])\in X\times \Pic^{\delta}(Y)|\ H^0({\cal M} \otimes {\cal E}^x)\ne
0\}\subset X\times \Pic^{\delta}(Y)  \ .
$$
$W$ is a non-empty   Zariski open set in $X$, whereas $Z$ is a closed analytic subset of $X\times
\Pic^{\delta}(Y)$, which obviously dominates
$X$. Choose any irreducible component $Z^0$ of $Z$ which still dominates $X$  and a    non-empty 
Zariski open set 
$V\subset X$ over which the fibres $Z^0_x$   have minimal
dimension. Choosing in a suitable way the multiplicities of the fibre components, one gets a holomorphic map
$V\ni x\mapsto Z_x^0$ with values in  the Chow variety of
$\Pic^{\delta}(Y)$ (\cite{Ba2}, p. 27). This
map has a meromorphic extension\footnote{The meromorphic extension theorem used here has been proved
in full generality in \cite{Ba1}; the proof is much easier when the target manifold is
K\"ahlerian (as is $\Pic^ \delta(Y)$ in our case). I am indebted to D. Barlet for explaining me
these results.}, which must be constant, because  
$a(X)=0$ and
this Chow variety  is algebraic.  Let
$z_0\subset \Pic^{\delta}(Y)$ be   this constant, and choose $[{\cal M}_0]\in z_0$. 
For $x\in V$ one has $Z_x^0=z_0$, whereas for $x\in W$ one has $Z_x^0\subset Z_x= z({\cal
E}^x)$. Therefore
\begin{equation}\label{zz}
[{\cal M}_0]\in z({\cal E}^x)\ \ \forall x\in V\cap W\ .
\end{equation}

Finally, let $U\subset V\cap W$ be the  Zarisky open  set of points $x\in V\cap W$ for
which 
$h^0({\cal M}_0\otimes {\cal E}^x)$ is minimal. By   (\ref{P})  and (\ref{zz}), this
minimal value    is either 1 or 2.   Set ${\cal T}:=(p_X)_*( p_Y^*({\cal M}_0)\otimes{\cal
E})$.  Using Grauert's local triviality and base change theorems (see
\cite{BHPV} p. 33), we see that
${\cal T}$ is locally free on $U$ and its fibre ${\cal T}(x)$ at a point $x\in U$ is $H^0({\cal
M}_0\otimes {\cal E}^x)$. Consider the natural morphism
$$e:p_X^*({\cal T})=p_X^*\left[(p_X)_*(p_Y^*({\cal M}_0)\otimes{\cal E})\right] \map 
p_Y^*({\cal M}_0)\otimes {\cal E} \ .
$$

For a point $(y,x)\in Y\times U$ the induced linear map  $e(y,x)$ between the corresponding fibres
is just the evaluation morphism   
$$H^0({\cal M}_0 \otimes {\cal E}^x) 
\map {\cal M}_0(y)\otimes{\cal E}^x(y)={\cal M}_0(y)\otimes{\cal E}_y(x)\ ,$$
which is injective  (by (\ref{P}) and (\ref{zz})). It suffices
now to consider the   morphisms 
$$e_y:{\cal T}  \map {\cal M}_0 (y)\otimes_\C{\cal E}_y \simeq{\cal E}_y  $$
obtained by
restricting
$e$ to
$\{y\}\times X$ and to take into account that the sheaf $p_X^*({\cal T})$ is locally free on $Y 
\times U$.
\qed
%

\begin{co}\label{type} In the hypothesis and with the notations of Theorem
\ref{constant} the following  holds: Either
${\cal E}_y$ is filtrable for every $y\in Y$, or ${\cal E}_y$ is non-filtrable for
every $y\in Y$.
\end{co}
\pf When $\rk({\cal T})=1$, all ${\cal E}_y$ will be filtrable.
When  $\rk({\cal T})=2$, the bundles ${\cal E}_y$ will be either all filtrable (when
${\cal T}$ is filtrable), or all non-filtrable (when ${\cal T}$ is
non-filtrable). 
\qed
\vspace{1mm}\\ 
We can now complete   the proof of our main theorem:
\\ \\
\pf (of Theorem \ref{main})
Let $X$ be a class  $VII$ surface such that none of the rational classes  $\pm c_1^\Q({\cal
K}),0,2c_1^\Q({\cal K})$ is  represented by a curve. Theorem \ref{curve} would yield  a compact
complex curve
$Y$ and a holomorphic morphism
$\fg:Y\to {\cal M}^{\rm st}$ taking both filtrable and non-filtrable values. But a minimal class $VII$
surface with positive $b_2$  has vanishing algebraic dimension. Therefore, by Lemma
\ref{family} and Corollary \ref{type}, such a morphism cannot exist.
\qed
\vspace{1mm}\\ 
{\bf Acknowledgments:} I am  grateful to Nicholas Buchdahl 	and  Matei Toma
for their encouragements and for the    useful discussions we had on
the preliminary version of the article.  
I   thank the two  referees for their   careful and valuable  suggestions. Thanks
to their comments, the presentation has been improved, several   errors and misprints
have been fixed and the last section has been substantially simplified.   I also thank G. Dloussky, K.
Oeljeklaus, F. Campana and D. Barlet for their interest in my work and useful dicussions on the
subject.
\vspace{3mm} 
{\small
Author's address: \vspace{2mm}\\
Andrei Teleman, LATP, CMI,   Universit\'e de Provence,  39  Rue F.
Joliot-Curie, 13453 Marseille Cedex 13, France,  e-mail:
teleman@cmi.univ-mrs.fr. }

\end{document}